\documentclass[11pt,a4paper]{article}
\usepackage{amsmath, amsfonts, amssymb,}
\usepackage{a4wide}
\usepackage{parskip}
\usepackage{enumitem}
\DeclareMathOperator{\spn}{span}
\def\proof{\noindent{\textbf{Proof. }}}
\def\QED{\hfill {$\square$}\goodbreak \medskip}

\newtheorem{Theorem}{Theorem}[section]
\newtheorem{Lemma}[Theorem]{Lemma}
\newtheorem{Proposition}[Theorem]{Proposition}

\newtheorem{Remark}[Theorem]{Remark}
\newtheorem{Definition}[Theorem]{Definition}

\numberwithin{equation}{section}

\begin{document}
	\date{}
	{\vspace{0.01in}
	\title{Three solutions for a fractional elliptic problem with asymmetric critical Choquard nonlinearity}
	\author{{\bf S. Rawat\footnote{email: {\tt sushmita.rawat1994@gmail.com}}  \;and  \bf K. Sreenadh\footnote{	e-mail: {\tt sreenadh@maths.iitd.ac.in}}} \\ Department of Mathematics,\\ Indian Institute of Technology Delhi,\\Hauz Khaz, New Delhi-110016, India. }
	\maketitle
	\begin{abstract}
	In this paper we study the existence and multiplicity of weak solutions for the following asymmetric nonlinear Choquard problem on fractional Laplacian:
	\begin{equation*}
	 \begin{array}{rl}
			 (-\Delta)^s u &= \displaystyle-\lambda|u|^{q-2}u + au + b\left( \int\limits_{\Omega}  \frac{(u^{+}(y))^{2^{*}_{\mu ,s}}}{|x-y|^ \mu}\, dy\right)  (u^{+})^{2^{*}_{\mu ,s}-2}u \quad\text{in} \; \Omega,\\
			 u &= 0\quad \text{in} \; \mathbb{R}^{N}\backslash\Omega,
		\end{array}
	\end{equation*}
where $\Omega$ is open bounded domain of $\mathbb{R}^{N}$ with $C^2$ boundary, $N > 2s$ and $s \in (0,1)$. Here $(-\Delta)^s$ is the fractional Laplace operator, $\lambda > 0$ is a real parameter, $q \in (1, 2)$, $a > 0$ and $b> 0$ are given constants, and $2^{*}_{\mu ,s} = \frac{2N-\mu}{N-2s}$ is the critical exponent in the sense of Hardy-Littlewood-Sobolev inequality and the notation $u^{+} = \max \{u, 0\}$. We prove that the above problem has at least three nontrivial solutions using the Mountain pass Lemma and Linking theorem.\medskip


\medskip


\end{abstract}
	
	\section{Introduction}
The purpose of this article is to study a doubly non-local Dirichlet problem with sub-linear and linear parts. We consider the problem
\begin{equation*}    (P_\lambda)\; \left\{  \begin{array}{rl}
		 (-\Delta)^s u &= \displaystyle-\lambda |u|^{q-2}u + au + b\left( \int\limits_{\Omega} \frac{(u^+(y))^{2^{*}_{\mu ,s}}}{|x-y|^ \mu}\,dy \right)  (u^+)^{2^{*}_{\mu ,s}-1}\quad \text{in} \; \Omega,\\
		 u &= 0\quad \text{in} \; \mathbb{R}^{N}\backslash\Omega,
	\end{array} \right.
\end{equation*}
where $\Omega$ is open bounded domain of $\mathbb{R}^{N}$ having $C^2$ boundary, $N > 2s$ with $s \in (0,1)$, $\lambda > 0$ is a real parameter, $q \in (1, 2)$, $a, b > 0$, where $2^{*}_{\mu ,s}= \frac{2N-\mu}{N-2s}$. The essential condition that we
assume is $\lambda_{k,s} < a < \lambda_{k+1,s}$, where $\{\lambda_{k,s}\}_{k\geq1}$ denote the sequence of eigenvalues of $(-\Delta)^s$  and $0<\mu < \min\{N,4s\}$. The notation $u^+ = \max\{u, 0\}\; \text{and}\; u^- = \min\{-u, 0\} $. Here, $(-\Delta)^s$ is fractional Laplace operator defined as,
\begin{equation*}
	(-\Delta)^su(x) := C(N, s) \lim \limits_{\epsilon \to 0} \int\limits_{\mathbb{R}^N\backslash B_{\epsilon}(x)}\dfrac{u(x)-u(y)}{|x-y|^{N+2s}}\,dy, \quad x\in \mathbb{R}^N,
\end{equation*}
where $C(N,s)$ is the normalization constant and $B_{\epsilon}(x)$ denotes the open ball of radius $\epsilon$ centred at $x$
\begin{equation*}
	C(N,s) = \left(\int\limits_{\mathbb{R}^{N}}\frac{1-cos(x_1)}{|x|^{N+2s}}\,dx\right) ^{-1}.
\end{equation*}
Due to the vast applications of non local operators, including it being the infinitesimal generators of the $L\acute{e}vy$ stable diffusion process \cite{Levy}, a lot of research has been done on this operator. For a comprehensive guide to the fractional operator we refer the interested readers to \cite{Hitchhiker's, fractionalS}. Servadei-Valdinoci \cite{Serv} investigated the non-local fractional
counterpart of the Laplace equation involving critical non-linearities, which was an extension of the famous work of Brezis and Nirenberg \cite{BN}. 
 Author in \cite{ser} extended the result obtained by Capozzi, Fortunato and Palmieri\cite{cap} and by Gazzola and Ruf \cite{gazzola}, to the non-local fractional counterpart of the Laplace equation with critical non-linearities. Objective of their work was to complete the existence result obtained in \cite[Theorem 4]{Serv} for any $\lambda > 0$ different from the eigenvalues of $(-\Delta)^s$.\par
In the local setting various authors have studied semi-linear elliptic problem of the type, for $1< q <2$
\begin{equation*}
	-\Delta u = -\lambda|u|^{q-2}u + g(u) \quad \text{in} \,\, \Omega
\end{equation*} 
where $g$ is asymmetric and asymptotically linear. Here asymmetric means that
$g$ satisfies an Ambrosetti–Prodi type condition  $\left( \text{i.e.}\,\, \displaystyle g_- := \lim\limits_{t \to -\infty } \frac{g(t)}{t} < \lambda_{k} < g_+ := \lim\limits_{t \to \infty } \frac{g(t)}{t}\right) $. Paiva-Massa \cite{paivamassa} studied this problem when $g$ is asymmetric and superlinear at $+ \infty$, $g_+ = \infty $ and assumed that $\displaystyle \frac{g(t)}{t}$ crosses an eigenvalue of the Laplacian when $t$ varies from $0$ to $-\infty $ (i.e. $g'(0) < \lambda_{k} < g_-)$. De Figueiredo-Yang \cite{fig} initiated the study of problem with critical exponent and asymmetric non-linearities to explore Ambrosetti-Prodi type problem for critical growth. Paiva-Presoto \cite{paiva} studied the critical problem for the case $g(u) = au + b(u^+)^{p-1}$ where $a \in \mathbb{R}$, $\lambda_{k} < a <\lambda_{k+1}$ and $b > 0$. Here $\{\lambda_{k}\}_{k\geq 1}$ denote the sequence of eigenvalues of $\left( -\Delta, H^1_0(\Omega)\right) $.\par
On the other hand one of the first applications of Choquard type non-linearity was given by Pekar in the framework of quantum theory \cite{pekar} and Lieb \cite{choqlieb} used it in approximation of Hartree-Fock theory. An extensive study has been done in the existence and uniqueness results due to its vast applications in physical models. For detailed understanding one can refer \cite{Moroz4, Moroz3, Moroz2, Moroz1} and references therein. In the local case $(s=1)$, Gao and Yang \cite{M.yang} analyzed Brezis-Nirenberg type existence results for Choquard type critical nonlinearity. In \cite{tuhinaFC}, authors studied these results in the following  non-local problem:
\begin{equation*}\label{ts}
	(-\Delta)^s u = \lambda u + \left( \int\limits_{\Omega}  \frac{|u(y)|^{2^{*}_{\mu, s}}}{|x-y|^ \mu}\, dy\right)  |u|^{2^{*}_{\mu,s}-2}u \quad\text{in} \; \Omega,\; u = 0 \; \text{on} \;\mathbb{R}^{N} \backslash \Omega, 
\end{equation*}
where $\Omega$ is a bounded domain in $\mathbb{R}^n$ with Lipschitz boundary, $\lambda$ is a real parameter, ${s\in(0, 1)}$, $N > 2s$, $ 0 < \mu < N$.
Gao and Yang \cite{gaoyang} established an existence and multiplicity result for the critical Choquard equation under the perturbation that is both sub-linear and sup-linear subcritical terms. 
They proved the existence of at least two positive solutions for an admissible small range of $\lambda$. Using sub- and super-solution techniques they obtained first solution which is local minimum for the functional. The second solution is obtained using Mountain-Pass lemma by studying the compactness of Palais-Smale sequences.\par
The motivation behind our paper is the inspiring work of  Miyagaki et al. \cite{miyagaki} where the authors have studied the following non-local problem consisting of sub-linear, linear and (critical) super-linear term,
\begin{equation*}
\begin{array}{cc}
	\displaystyle (-\Delta)^s u = -\lambda|u|^{q-2}u + au + b(u^{+})^{2^{*}_{s}-1} \;\;\text{in} \; \Omega,\;\;
	u = 0\quad \text{in} \; \mathbb{R}^{N}\backslash\Omega,
\end{array}
\end{equation*}
where $\Omega$ is open bounded domain of $\mathbb{R}^{N}$ having $C^2$ boundary, $N > 2s$ with $s \in (0,1)$, $\lambda > 0$ is a real parameter, $q \in (1, 2)$, $a, b > 0$, where $2^{*}_{s}= \frac{2N}{N-2s}$ and $\lambda_{k,s} < a < \lambda_{k+1,s}$. They have proved that for sufficiently small $\lambda >0$ the problem has atleast three non-trivial distinct solutions. First two solutions of opposite signs are obtained using Mountain-Pass theorem and Linking theorem is used to obtain the third solution.\par
To the best of our knowledge there is no work on doubly non-local case with sub-linear and linear term and even for the corresponding problem in the local case. In this article we study the multiplicity result of $(P_\lambda)$. Due to the lack of compactness we will be proving Palais-Smale condition is satisfied below some suitable critical level. For the non-trivial non-negative solution we will satisfy the Mountain-Pass geometry for some suitable choice of $\lambda > 0$. Whereas for the non-trivial non-positive solution we will satisfy the Mountain-Pass geometry for all $\lambda > 0$. Next we use minimax result of Generalized Mountain-Pass theorem to obtain third non-trivial solution. With this introduction we state our main result.
\begin{Theorem}\label{thm1}
	Let $1 <q< 2$, $\lambda_{k,s} < a < \lambda_{k+1,s}$, $b > 0$, $0 < s < 1$ and $N \geq 4s$. Then
	there exists $\lambda > 0$ such that if $0 < \lambda < \overline{\lambda}$, problem $(P_\lambda)$ possesses at least three nontrivial solutions
	with one non-negative and one non-positive. The same conclusion holds for the case $2s<N< 4s$ if the
	hypothesis $\lambda_{k,s} < a < \lambda_{k+1,s}$ is met with an integer $k$ large enough.
\end{Theorem}
\begin{Remark}
Our results will hold for the following local case problem with suitable modifications of our analysis and method
\begin{equation*}
	\begin{array}{rl}
		-\Delta u &= \displaystyle-\lambda|u|^{q-2}u + au + b\left( \int\limits_{\Omega}  \frac{(u^{+}(y))^{2^{*}_{\mu }}}{|x-y|^ \mu}\, dy\right)  (u^{+})^{2^{*}_{\mu }-2}u \quad\text{in} \; \Omega,\\
		u &= 0\quad \text{on} \; \partial\Omega,
	\end{array}
\end{equation*}
where $\Omega$ is open bounded domain of $\mathbb{R}^{N}$ having $C^2$ boundary, $N\geq 3$, $\lambda > 0$ is a real parameter, $q \in (1, 2)$, $a, b > 0$, where $2^{*}_{\mu }= \frac{2N-\mu}{N-2}$. Assuming that $\lambda_{k} < a < \lambda_{k+1}$, where $\{\lambda_{k}\}_{k\geq1}$ denote the sequence of eigenvalues of $(-\Delta, H^1_0(\Omega))$ and $0<\mu < \min\{N,4\}$.
\end{Remark}
The paper is organized as follows: In section 2 we present some preliminaries on function spaces required for variational settings. In section 3 we prove the result regarding Palais-Smale condition. In section 4 we infer the existence of two opposite sign solutions. In section 5 we obtain results required for the Linking theorem. In section 6 we prove Theorem \ref{thm1} and show the existence of the third solution by estimating the minimax value obtained by Linking theorem.

    \section{Preliminaries}
  We recall some definitions of function spaces and results that will be required in later sections.
  Consider the functional space $H^s(\mathbb{R}^N)$ as the usual fractional Sobolev space  defined as
  \begin{equation*}
  	H^s(\mathbb{R}^N) = \left\lbrace u\in L^2(\mathbb{R}^N): \left[ u\right]_s < \infty \right\rbrace,
  \end{equation*}
where $\left[ u\right]_s := \left( \displaystyle\int\limits_{\mathbb{R}^{2N}} \frac{|u(x)-u(y)|^2}{|x-y|^{N+2s}}\,dxdy\right)^{\frac{1}{2}}$ is the Gagliardo seminorm of measurable function $u$
  and the space is endowed with the norm 
  \begin{equation}\label{eq2.1}
  	\|u\|_{H^s(\mathbb{R}^N)} = \left(\|u\|^{2}_{L^2(\mathbb{R}^N))} + \left[ u\right]_s^2 \right)^{\frac{1}{2}}.
  \end{equation}
  We define the functional space associated to this problem as
  \begin{equation*}
  	X_0 = \left\lbrace u\in H^s(\mathbb{R}^N): u = 0\; a.e.\; in \; \mathbb{R}^N\backslash \Omega  \right\rbrace
  \end{equation*} 
which is a closed subspace of $H^s(\mathbb{R}^N)$.
  Then, it can be shown that $X_0$ is a Hilbert space with the inner product 
  \begin{equation*}
  	\langle u, v \rangle = \int\limits_{\mathbb{R}^{2N}}\dfrac{(u(x)-u(y))(v(x)-v(y))}{|x-y|^{N+2s}}\,dxdy,
  \end{equation*}
  for $u,v \in X_0$ and thus the corresponding norm,
  \begin{equation*}
  	\|u\|_{X_0} = \|u\| = \displaystyle \left( \int\limits_{\mathbb{R}^{2N}} \frac{|u(x)-u(y)|^2}{|x-y|^{N+2s}}\,dxdy \right)^{\frac{1}{2}}.
  \end{equation*} It can be shown that this is equivalent to \eqref{eq2.1} on $X_0$.
 We refer to \cite[Proposition 3.6]{Hitchhiker's} for the equality
\begin{equation*}
	[u]_s^2 = 2C(N,s)^{-1}\int\limits_{\mathbb{R}^{N}}|(-\Delta)^{\frac{s}{2}}u|^2\,dx,
\end{equation*} 
for all $u \in H^s(\mathbb{R^N})$.
Thus, for all $u, v \in X_0$ there holds
\begin{equation*}
	2C(N, s)^{-1}\int\limits_{\mathbb{R}^{N}}u(x)(-\Delta)^sv(x)\,dx = \int\limits_{\mathbb{R}^{2N}}\dfrac{(u(x)-u(y))(v(x)-v(y))}{|x-y|^{N+2s}}\,dxdy.
\end{equation*}
In particular, it follows that the linear operator $(-\Delta)^s$ is self-adjoint on $X_0(\Omega$). $\{\lambda_{j,s}\}_{j\geq1}$ denote the sequence of eigenvalues of $(-\Delta)^s$ on the space $X_0$ satisfying
\begin{equation*}
0 < \lambda_{1,s} <\lambda_{2,s} \leq \lambda_{3,s} \leq \cdots \leq \lambda_{k,s} \leq \lambda_{k+1,s} \cdots,\, \lambda_{k,s} \to \infty \;\;\text{as}\;\; j\to \infty,
\end{equation*}
where each eigenvalue $\lambda_{k,s}$ is repeated according to its multiplicity and $\phi_{k,s} \in C^{0, \sigma}$ is the eigen function associated to the eigenvalue $\lambda_{k,s}$, with some $\sigma \in (0, 1)$. For the rest of the paper we fix a sequence $\{\phi_{k,s}\}$ of eigen functions forming an orthonormal basis in $L^2(\Omega)$ and an orthogonal basis in $X_0(\Omega)$. Furthermore, The smallest eigenvalue $\lambda_1$ is simple and isolated in the spectrum and we assume $\phi_{1,s} > 0$ in $\Omega$.
  \begin{Proposition}
  	\textbf{Hardy-Littlewood-Sobolev inequality}: Let $t$, $r > 1$ and $0 < \mu < N$ with $\frac{1}{t} + \frac{\mu}{N} + \frac{1}{r} = 2$, $f \in L^t(\mathbb{R}^N)$ and $h \in L^r(\mathbb{R}^N)$. Then there exists a sharp constant $C(t, r, \mu, N)$ independent of $f$, $h$ such that
  	\begin{equation*}
  		\iint\limits_{\mathbb{R}^{2N}} \dfrac{f(x)h(y)}{|x-y|^ \mu}\,dxdy \leq C(t, r, \mu, N)\|f\|_{L^t(\mathbb{R}^N)}\|h\|_{L^r(\mathbb{R}^N)}.
  	\end{equation*}
  	
  \end{Proposition}
We are looking for nontrivial solutions of the following problem $(P_\lambda)$ 
and its weak solution which is defined as, $u \in X_0$ such that
\begin{equation*}\label{eq2.3}
	\begin{aligned}
		\displaystyle \frac{C(N,s)}{2}\int\limits_{\mathbb{R}^{2N}}\frac{(u(x)-u(y))(\phi(x)-\phi(y))}{|x-y|^{N+2s}}\,dxdy = &-\lambda \int\limits_{\Omega}|u(x)|^{q-2}u(x)\phi(x)dx + a\int\limits_{\Omega}u\phi(x)dx\\ &+b\int\limits_{\Omega}\int\limits_{\Omega}\frac{(u^{+}(y))^{2^{*}_{\mu ,s}}(u^+(x))^{2^{*}_{\mu ,s}-1}\phi(x)}{|x-y|^ \mu}\,dxdy
	\end{aligned}
\end{equation*}
for any $\phi \in X_0$. 
 The energy functional associated with the problem $(P_\lambda)$ is ${J_{\lambda } : X_0(\Omega) \rightarrow \mathbb{R}}$\; defined as,
\begin{equation}\label{eq2.4}
	\begin{aligned}
		\displaystyle J_\lambda (u)= &\frac{C(N,s)}{4} \int\limits_{\mathbb{R}^{2N}} \frac{|u(x)-u(y)|^2}{|x-y|^{N+2s}}\,dxdy +\frac{\lambda}{q} \int\limits_{\Omega} |u(x)|^{q} \,dx -\frac{a}{2}\int\limits_{\Omega} |u(x)|^{2}\,dx\\
		&-\frac{b}{2\cdot{2^*_{\mu ,s} }}\int\limits_{\Omega}\int\limits_{\Omega} \frac{(u^{+}(y))^{2^{*}_{\mu ,s}}(u^{+}(x))^{2^{*}_{\mu ,s}}}{|x-y|^ \mu}\,dxdy,
	\end{aligned}
\end{equation}
which is $C^1$ functional. For $\phi \in X_0(\Omega)$ we define the Gateaux derivative as,
\begin{equation}\label{eq2.5}
	\begin{aligned}
		\displaystyle \langle J'_\lambda (u), \phi\rangle= &\frac{C(N,s)}{2} \int\limits_{\mathbb{R}^{2N}} \frac{(u(x)-u(y))(\phi(x)-\phi(y))}{|x-y|^{N+2s}}\,dxdy +\lambda \int\limits_{\Omega} |u(x)|^{q-2}u(x)\phi(x) \,dx \\&-a\int\limits_{\Omega} u(x)\phi(x)\,dx
		-b\int\limits_{\Omega}\int\limits_{\Omega} \frac{(u^{+}(y))^{2^{*}_{\mu ,s}}(u^{+}(x))^{2^{*}_{\mu ,s-1}}\phi(x)}{|x-y|^ \mu}\,dxdy.
	\end{aligned}
\end{equation}
From the embedding results, we conclude that $X_0$ is continuously embedded in $L^{p}(\Omega)$ when $1 \leq p \leq 2_s^{*}$. Also the embedding is compact for $1 \leq p < 2_s^{*}$, but not for the case $p = 2_s^{*}$. We define
the best constant for the embedding $X_0$ into $L^{2_s^{*}}(\mathbb{R}^N)$ as, 
\begin{equation*}\label{eq2.7}
	S_s =\inf\limits_{u\in X_0\backslash \{0\}}\left\lbrace\int\limits_{\mathbb{R}^{2N}} \frac{|u(x)-u(y)|^2}{|x-y|^{N+2s}}\,dxdy:  \int\limits_{\mathbb{R^N}}|u|^{2_s^{*}} = 1 \right\rbrace.
\end{equation*}
Consequently, we define
\begin{equation}\label{eq2.8}
	S_s^H = \inf\limits_{u\in X_0\backslash \{0\}}\left\lbrace\int\limits_{\mathbb{R}^{2N}} \frac{|u(x)-u(y)|^2}{|x-y|^{N+2s}}\,dxdy:  \int\limits_{\mathbb{R}^{2N}} \dfrac{|u(x)|^{2^{*}_{\mu ,s}}|u(y)|^{2^{*}_{\mu ,s}}}{|x-y|^ \mu}\,dxdy = 1\right\rbrace. 
\end{equation}
\begin{Lemma} \cite{Tuhina}
	The constant $S_s^H$ is achieved by u if and only if u is of the form\\
	$C\left( \frac{t}{t^2 + |x-x_0|^2}\right) ^\frac{N-2s}{2}$, $x \in \mathbb{R}^N$, 
	for some $x_0 \in \mathbb{R}^N, C \text{and}\; t > 0.$ Moreover,
	$S_s^H = \frac{S_s}{{C(N, \mu)}^\frac{1}{2_{\mu, s}^{*}}}$.
\end{Lemma}

Consider the family of functions ${U_\epsilon}$, where $U_\epsilon$ is defined as
\begin{equation} \label{minmimizer}
	U_\epsilon = \epsilon^{-\frac{N-2s}{2}}u^{*}\left( \frac{x}{\epsilon}\right), \; x\in \mathbb{R}^{N}, \epsilon > 0,
\end{equation}
\begin{equation*}
	u^{*}(x) = \overline{u}\left( \frac{x}{{S_s}^\frac{1}{2s}}\right) , \; \overline{u}(x) = \frac{\tilde{u}(x)}{\|\tilde u\|_{L^{2_s^{*}}(\mathbb{R}^N)}} \; \text{and}\; \tilde{u}(x) = \alpha(\beta^2 + |x|^2)^{-\frac{N-2s}{2}},
\end{equation*}
with $\alpha \in \mathbb{R}\backslash\{0\}$ and $\beta > 0$ are fixed constants.
Then for each $\epsilon > 0, \; U_\epsilon $ satisfies
\begin{equation*}
	(-\Delta)^s u = |u|^{2_s^{*}-2}u \quad in \; \mathbb{R}^N,
\end{equation*}
and the equality,
\begin{equation*}
	\int\limits_{\mathbb{R}^{2N}} \frac{|U_\epsilon(x)-U_\epsilon(y)|^2}{|x-y|^{N+2s}}\,dxdy = \int\limits_{\mathbb{R^N}}|U_\epsilon|^{2_s^{*}} = {S_s}^\frac{N}{2s}. 
\end{equation*}
We recall 
\begin{Definition} Let $Y$ be a real Banach space. Then
	\begin{enumerate}
		\item A sequence  $(u_n)\subset Y $ is called a Palais-Smale sequence for the functional $I \in C^1(Y,\mathbb{R})$ at the level $r \in \mathbb{R}$ if $I(u_n) \to r$ as $n \to \infty$, and $|\langle I'(u_n), \phi\rangle|_{Y} \leq \epsilon_n \|\phi\|_{Y}$ for all $\phi \in Y$, where $\epsilon_n \to 0$ as $n \to \infty$.
		\item $I$ is said to satisfy Palais-Smale condition at level $r$ if every Palais-Smale sequence at level $r$ has a convergent subsequence. 
	\end{enumerate}
	
\end{Definition}
  \begin{Theorem}\label{thm3}
  		Let $Y$ be a real Banach space satisfying $Y = V\oplus W$ with $V$ finite dimensional. Suppose that the functional $I \in C^1(Y,\mathbb{R})$ fulfills the conditions:
  		\begin{enumerate}
  	\item there are constants $\rho, \alpha > 0$ such that $I_{\partial B_\rho \cap W} \geq \alpha$, where $B_\rho$ denotes the open ball centered
  	at zero and of radius $\rho$.
  	\item there are constants $R_1, R_2 > \rho$, $\beta < \alpha$ and a nonzero vector $e  \in W$ such that $I_{\partial Q} \leq \beta$, with
  	$Q = \{v + re : v \in B_{R_1} \cap V , 0 < r < R_2\}$.
  	 \end{enumerate}
  	Then $I$ possesses a Palais-Smale sequence at the level $c \geq \alpha$ where,
  	$c = \inf\limits_{h \in \Gamma}\max\limits_{u \in Q} I(h(u))$,
  	and $\Gamma = \{h \in C(\overline{Q}, Y ) : h = \text{id}\; \text{on}\; \partial Q\}$. \QED
  \end{Theorem}
Throughout the paper we will use the notation,
\begin{equation*}
	\|u\|_{0}^{2\cdot2^{*}_{\mu ,s}} := \int\limits_{\Omega}\int\limits_{\Omega} \frac{|u(x)|^{2^{*}_{\mu ,s}}|u(y)|^{2^{*}_{\mu ,s}}}{|x-y|^ \mu}\,dxdy.
\end{equation*}
  
  \section{Palais-Smale condition}
  In this section  we study the compactness of Palais-Smale sequences of the functional $J_\lambda$. 
\begin{Lemma}\label{lemma3.1}
	Let $\lambda > 0$, $1< q <2$, $\lambda_{k,s} < a < \lambda_{k+1,s}$, $b>0$. Then the functional $J_\lambda$ satisfies Palais-Smale condition at level $c < \left( \dfrac{2^{*}_{\mu,s}-1}{2\cdot2^{*}_{\mu,s}} \right)\left[\dfrac{C(N,s)S_s^H}{2}\right]^{\frac{2^{*}_{\mu,s}}{2^{*}_{\mu,s}-1}}b^{\frac{-1}{2^{*}_{\mu,s}-1}} $.
\end{Lemma}
	\proof Let $\langle u_n \rangle$ be palais-smale sequence in $X_0(\Omega)$ for the functional $J_\lambda$. Then there exists a positive constant $C$ such that
	 $|J_\lambda(u_n)| \leq C$ and $|\langle J'_\lambda(u_n), \phi\rangle| \leq \epsilon_n \lVert\phi\rVert$ for all $\phi \in X_0(\Omega)$ and $\epsilon_n \to 0$ as $n \to \infty$. We have to show that $\langle u_n \rangle$ has a convergent subsequence. First we will show that $\lVert\langle u_n \rangle\rVert$ is bounded. Let if possible $\lVert\langle u_n \rangle\rVert \to \infty$ as $n \to \infty$. By \eqref{eq2.4} and \eqref{eq2.5} we get,
	 \begin{equation*}
	 	\begin{aligned}
	 		J_\lambda(u_n)-\dfrac{\langle J'_\lambda(u_n), u_n\rangle}{2} &= \lambda\left[\dfrac{1}{q}-\dfrac{1}{2} \right]\int\limits_{\Omega}|u_n(x)|^{q}\,dx + b\left[ \dfrac{1}{2}-\dfrac{1}{2\cdot2^{*}_\mu,s}\right]\lVert u^{+}_n \rVert_0^{2^{*}_{\mu,s}}\\
	 		&\geq \lambda\left[\dfrac{1}{q}-\dfrac{1}{2} \right]\int\limits_{\Omega}|u_n(x)|^{q}\,dx.
	 	\end{aligned}
	 \end{equation*}
	 As $1< q < 2$ and by Palais-Smale condition we imply, 
	 \begin{equation}\label{eq3.1}
	 	0 \leq \lambda\left[\dfrac{1}{q}-\dfrac{1}{2} \right]\int\limits_{\Omega}|u_n(x)|^{q}\,dx \leq C+ \epsilon_n\lVert u_n\rVert.
	 \end{equation}
	 Set $v_n:= \dfrac{u_n}{\lVert u_n\rVert}$, and since $\lVert v_n \rVert = 1$, it is a bounded sequence in $X_0(\Omega)$. Thus upto  a subsequence, there exists a function $v \in X_0(\Omega)$ such that as $n \to \infty$, we have $v_n \rightharpoonup v$ weakly in $X_0(\Omega)$, $v_n \rightarrow v$ strongly in $ L^r(\Omega)$ where $r \in [1, 2^{*}_{s})$ and $v_n \rightarrow v$ a.e in $\Omega$.
	 Now dividing \eqref{eq3.1} by $\lVert u_n\rVert$ and taking $n \to \infty$ we get,
	 \begin{equation*}
	 	\lim\limits_{n \to \infty}\dfrac{\int\limits_{\Omega}|v_n(x)|^{q}\,dx}{{\lVert u_n \rVert}^{1-q}} = 0.
	 \end{equation*}
	 i.e \begin{equation*}
	 \lim\limits_{n \to \infty}{\lVert u_n\rVert}^{q-1}	\int\limits_{\Omega}|v|^q\,dx = 0.
	 \end{equation*}
	 If $v \neq 0$, then that gives us a contradiction. Thus, we assume $v=0$
	 \begin{equation*}
	 	v_n \rightharpoonup 0\quad \text{in}\;\; X_0(\Omega).
	 \end{equation*}
 Next we consider
 \begin{equation*}
 \begin{aligned}
 	J_\lambda(u_n)-\dfrac{\langle J'_\lambda(u_n), u_n\rangle}{2\cdot2^*_{\mu,s}} =&\left[\frac{1}{2}-\frac{1}{2\cdot2^*_{\mu,s}} \right] \|u_n\|^2+ \lambda\left[\dfrac{1}{q}-\dfrac{1}{2\cdot2^*_{\mu,s}} \right]\int\limits_{\Omega}|u_n(x)|^{q}\,dx\\ 
 	&- a\left[ \dfrac{1}{2}-\dfrac{1}{2\cdot2^{*}_{\mu,s}}\right]\int\limits_{\Omega}|u_n(x)|^{2}\,dx.	
\end{aligned}
\end{equation*}
 Dividing by $\|u_n\|^2$ and taking $n \to \infty$ we get
 \begin{equation*}
 	\begin{aligned}
 		\left[\frac{1}{2}-\frac{1}{2\cdot2^*_{\mu,s}} \right] \left\lbrace 1-a\int\limits_{\Omega}|v(x)|^2\,dx \right\rbrace &= 0\\ 
 		\int\limits_{\Omega}|v(x)|^2\,dx  = \frac{1}{a} &> 0
 	\end{aligned}
 \end{equation*}
 which holds since $a >0$, but this contradicts the fact that $v = 0$. Hence our assumption is wrong and thus $\langle\lVert u_n \rVert\rangle$ is a bounded sequence.\\
	 Next we will show that $\langle u_n \rangle$ has a convergent subsequence. Since $X_0(\Omega)$ is a reflexive Banach space, thus upto  a subsequence, there exists a function $u \in X_0(\Omega)$ such that as $n \to \infty$, we have $u_n \rightharpoonup u$ weakly in $X_0(\Omega)$, $u_n \rightarrow u$ strongly in $ L^r(\Omega)$ where $r \in [1, 2^{*}_{s})$,$u_n \rightarrow u$ a.e in $\Omega$ and $\displaystyle{\int\limits_{\Omega} \frac{(u^+_n(y))^{2^{*}_{\mu ,s}}(u^+_n(x))^{2^{*}_{\mu ,s}-1}}{|x-y|^ \mu}\,dy \rightharpoonup \int\limits_{\Omega} \frac{(u^+(y))^{2^{*}_{\mu ,s}}(u^+(x))^{2^{*}_{\mu ,s}-1}}{|x-y|^ \mu}\,dy} $ weakly in $L^\frac{2N}{N+2s}$.\\
	 Note that $u$ is the weak solution of our Problem $(P_\lambda)$. So we imply,
	 \begin{equation}\label{eq3.2}
	 		J_\lambda(u)= J_\lambda(u)-\dfrac{\langle J'_\lambda(u), u\rangle}{2}\;\geq 0.
	 	 \end{equation}
	 From the (P.S) assumption we infer that,
	 \begin{equation}\label{eq3.3}
	 	c =\lim\limits_{n \to \infty}J_\lambda(u_n)= \lim\limits_{n \to \infty}\dfrac{C(N, s)}{4}{\| u_n\|}^2+ \dfrac{\lambda}{q}\|u_n\|^q_{L^q}-\dfrac{a}{2}\|u_n\|^2_{L^2}-\dfrac{b}{2\cdot2^{*}_{\mu,s}}\|u^+_n\|^{2\cdot2^{*}_{\mu,s}}_{0}.
	 \end{equation}
	 We set $v_n = u_n - u$, $v_n \rightharpoonup v $ in $X_{0}$. By Brezis-Lieb \cite{BrezLieb} and \cite{M.yang} we derive that for $n \to \infty$,
	 \begin{equation*}
	 \begin{split}
	 	\|u_n\|^2 &= \|v_n\|^2 + \|u\|^2 + o(1)\\
	 	\|u^+_n\|_0^{2\cdot2^{*}_{\mu, s}} &= \|v^+_n\|_0^{2\cdot2^{*}_{\mu, s}} + \|u^+\|_0^{2\cdot2^{*}_{\mu, s}} + o(1).
	 \end{split}
	 \end{equation*}
	 Using the above in \eqref{eq3.3} to get,
	 \begin{equation}\label{eq3.4}
	 	c = \lim\limits_{n \to \infty}J_\lambda(v_n)+ J_\lambda(u).
	 \end{equation}
	 By arguing as in \eqref{eq3.3} and using \eqref{eq3.2} we get,
	  \begin{equation*}
	  	\begin{split}
	  		0 &= \lim\limits_{n \to \infty}\langle J'_\lambda(u_n), u_n\rangle\\
	  		&= \lim\limits_{n \to \infty}\langle J'_\lambda(v_n), v_n\rangle\\
	  		&= \lim\limits_{n \to \infty}\frac{C(N,s)}{2}\|v_n\|^2-b\|v^+_n\|_0^{2\cdot2^{*}_{\mu,s}}.
	  	\end{split}
	  \end{equation*}
	  Fix a relabeled sequence such that, $$\lim\limits_{n \to \infty}\|v_n\|^2 = l$$
	  which implies, $$\lim\limits_{n \to \infty}\|v^+_n\|_0^{2} = \left( {\frac{C(N,s) l}{2b}}\right) ^\frac{1}{2^{*}_{\mu,s}}.$$
	  If suppose $l > 0$, then by Sobolev embedding \eqref{eq2.8} we have,
	  \begin{equation}\label{eq3.5}
	  	\begin{split}
	  		\|v_n\|^{2} &\geq S_s^H\|v^{+}_n \|_0^2\\
	  			l &\geq \left( S_s^H\right)^\frac{2^{*}_{\mu,s}}{2^{*}_{\mu,s}-1} \left( {\frac{C(N,s)}{2b}}\right) ^\frac{1}{2^{*}_{\mu,s}-1}.
	  		\end{split}
	  \end{equation}
  Combining \eqref{eq3.4}, \eqref{eq3.2} and \eqref{eq3.5} we get,
  \begin{equation*}
  	\begin{split}
  		c\geq \lim\limits_{n \to \infty}J_\lambda(v_n) &=\frac{C(N,s)}{4}l-\frac{b}{2\cdot2^{*}_{\mu,s}}\frac{C(N,s)l}{2b}\\
  		&\geq \left( \dfrac{2^{*}_{\mu,s}-1}{2\cdot2^{*}_{\mu,s}} \right)\left[\dfrac{C(N,s)S_s^H}{2}\right]^{\frac{2^{*}_{\mu,s}}{2^{*}_{\mu,s}-1}}b^{\frac{-1}{2^{*}_{\mu,s}-1}},
  	\end{split}
   \end{equation*}
which contradicts the assumption on $c$. Hence, $l = 0$ which further implies $u_n \rightarrow u$ strongly in $X_0(\Omega)$ upto a subsequence. Therefore, $J_\lambda$ satisfies Palais-Smale condition at any level \\$c < \left( \dfrac{2^{*}_{\mu,s}-1}{2\cdot2^{*}_{\mu,s}} \right)\left[\dfrac{C(N,s)S_s^H}{2}\right]^{\frac{2^{*}_{\mu,s}}{2^{*}_{\mu,s}-1}}b^{\frac{-1}{2^{*}_{\mu,s}-1}}.$\QED                                                      

\section{Constant sign solutions}
In this section we will find two non-trivial solution of opposite sign for our problem $(P_\lambda)$. It is derived with the help of $H^s$ v.s $C^0_s$ minimizer property proved by \cite{Divya do-reg}. Recall that with $\delta^s(x)$ := $dist(x, \partial\Omega)$, the space $C^0_s$ is defined as,
\begin{equation*}
	C^0_s(\overline{\Omega})= \{u \in C^0(\overline{\Omega}): \|\frac{u}{\delta^s}\|_{L^\infty} < \infty\},
\end{equation*}
endowed with the norm
\begin{equation*}
	\|u\|_{C^0_s} = \|\frac{u}{\delta^s}\|_{L^\infty}.
\end{equation*}
For the positive case the corresponding auxiliary problem will be
\begin{equation}\label{eq4.1}  \begin{aligned}
		\displaystyle (-\Delta)^s u &= -\lambda (u^+)^{q-1} + au^+ + b\left( \int\limits_{\Omega} \frac{(u^+(y))^{2^{*}_{\mu ,s}}}{|x-y|^ \mu}\,dy \right)  (u^+)^{2^{*}_{\mu ,s}-1}\quad \text{in} \; \Omega,\\
		u & = 0\quad \text{in} \; \mathbb{R}^{N}\backslash\Omega.
	\end{aligned}
\end{equation}
The energy functional associated with the problem \eqref{eq4.1} is ${J^{+}_{\lambda } : X_0(\Omega) \rightarrow \mathbb{R}}$ defined as,
\begin{equation*}\label{eq4.2}
	\begin{aligned}
		\displaystyle J^{+}_\lambda (u)= &\frac{C(N,s)}{4} \int\limits_{\mathbb{R}^{2N}} \frac{|u(x)-u(y)|^2}{|x-y|^{N+2s}}\,dxdy +\frac{\lambda}{q} \int\limits_{\Omega} (u^+(x))^{q} \,dx -\frac{a}{2}\int\limits_{\Omega} (u^+(x))^{2}\,dx\\
		&-\frac{b}{2\cdot{2^*_{\mu ,s} }}\int\limits_{\Omega}\int\limits_{\Omega} \frac{(u^{+}(y))^{2^{*}_{\mu ,s}}(u^{+}(x))^{2^{*}_{\mu ,s}}}{|x-y|^ \mu}\,dxdy.
	\end{aligned}
\end{equation*}
We observe that $J^+_\lambda$ is of $C^1$ class and its derivative is
\begin{equation}\label{eq4}
\begin{aligned}
	 \langle (J^{+})'_\lambda (u), \phi \rangle = &\frac{C(N,s)}{2} \int\limits_{\mathbb{R}^{2N}} \frac{(u(x)-u(y))(\phi(x)-\phi(y))}{|x-y|^{N+2s}}\,dxdy +\frac{\lambda}{q} \int\limits_{\Omega} (u^+(x))^{q-1}\phi(x) \,dx\\ &-\frac{a}{2}\int\limits_{\Omega} u^+(x)\phi(x)\,dx
	-\frac{b}{2\cdot{2^*_{\mu ,s} }}\int\limits_{\Omega}\int\limits_{\Omega} \frac{(u^{+}(y))^{2^{*}_{\mu ,s}}(u^{+}(x))^{2^{*}_{\mu ,s}-1}\phi(x)}{|x-y|^ \mu}\,dxdy.
\end{aligned}
\end{equation}
thus the critical point of $J^+_\lambda$ are the non negative weak solution for problem $(P_\lambda)$.
\begin{Theorem}\label{thm4.1}
	Let $\lambda \in (0, \lambda_0)$ for some $\lambda_0 > 0$. Then problem $(P_\lambda)$ has a positive solution in $X_0(\Omega)$.
\end{Theorem}
	\proof We will show there exist positive constants $\alpha$ and $\rho$ such that,
	\begin{equation}\label{eq4.3}
		J^+_\lambda(u) \geq \alpha \quad \text{for all}\;\; u\in X_{0}(\Omega),\; \|u\| = \rho.
	\end{equation}
	We construct a functional $I_\lambda(u)$ and claim $u =0$ is strict local minimizer i.e. $I_\lambda(u) \geq 0$ for all $\|u\| < \rho_1$ for some $\rho_1 >0$, , \begin{equation*}
		I_\lambda(u):= \frac{\lambda}{q} \int\limits_{\Omega} (u^+(x))^{q} \,dx -\frac{a}{2}\int\limits_{\Omega} (u^+(x))^{2}\,dx-\frac{b}{2\cdot{2^*_{\mu ,s} }}\int\limits_{\Omega}\int\limits_{\Omega} \frac{(u^{+}(y))^{2^{*}_{\mu ,s}}(u^{+}(x))^{2^{*}_{\mu ,s}}}{|x-y|^ \mu}\,dxdy.
	\end{equation*}
	By \cite{Divya do-reg} it suffices to prove that $0$ is the local minimizer for $I_\lambda$ on $X_0 \cap C^0_s(\overline{\Omega})$ i.e. $I_\lambda(u) \geq 0$ for all $\|u\|_{C_s^0} < \rho_2$ for some $\rho_2 >0$.
	 For any $u \in X_0 \cap C^0_s(\overline{\Omega})$ we imply that
	\begin{equation}\label{eq4.4}
    \begin{aligned}
    	\int\limits_{\Omega}(u^+(x))^{2}\,dx &= \int\limits_{\Omega}\left(\frac{u^+(x)}{\delta^s(x)}\right)^{2-q}{\left( \delta^s(x)\right)}^{2-q} (u^+(x))^{q}\,dx\\
    	&\leq C_1 \|u\|_{C^0_s}^{2-q}\int\limits_{\Omega}(u^+(x))^{q}\,dx.
    \end{aligned}
	\end{equation}
Since we are proving the claim for all $\|u\| < \rho_1$ for some $\rho_1 >0$, thus it will also work for sufficiently small ball i.e. for $\|u\|$ very small. By Sobolev-embedding we have $\|u\|_{L^{2^*_s}} \leq c \|u\|$ for some constant $c>0$.
Hence, we conclude by Hardy-Littlewood-Sobolev inequality,
\begin{equation}\label{eq4.5}
\begin{aligned}
	\int\limits_{\Omega}\int\limits_{\Omega} \frac{(u^{+}(y))^{2^{*}_{\mu ,s}}(u^{+}(x))^{2^{*}_{\mu ,s}}}{|x-y|^ \mu}\,dxdy &\leq c\|u^+\|^{2\cdot2^*_{\mu,s}}_{L^{2^*_s}}
	 \leq c\|u^+\|^{2^*_s}_{L^{2^*_s}}\\
	& \leq c \int\limits_{\Omega}\left(\frac{u^+(x)}{\delta^s(x)}\right)^{2^{*}_{s}-q}{\left( \delta^s(x)\right)}^{2^{*}_{s}-q} (u^+(x))^{q}\,dx\\
	& \leq C_2 \|u\|_{C^0_s}^{2^{*}_{s}-q}\int\limits_{\Omega}(u^+(x))^{q}\,dx,
\end{aligned}
\end{equation}
	where $C_1$ and $C_2$ are positive constants. From \eqref{eq4.4} and \eqref{eq4.5} we get,
	\begin{equation*}
		I_\lambda(u) \geq \left[ \frac{\lambda}{q} -\frac{a C_1}{2}\|u\|_{C^0_s}^{2-q} -\frac{b C_2}{2\cdot{2^*_{\mu ,s} }}\|u\|_{C^0_s}^{2^{*}_{s}-q}\right] \int\limits_{\Omega} (u^+(x))^{q} \,dx \geq 0,
	\end{equation*}
	whenever \begin{equation*}
		\frac{\lambda}{q} \geq \frac{a C_1}{2}\|u\|_{C^0_s}^{2-q} + \frac{b C_2}{2\cdot{2^*_{\mu ,s} }}\|u\|_{C^0_s}^{2^{*}_{s}-q}.
	\end{equation*}
	Since \begin{equation*}
			\displaystyle J^{+}_\lambda (u)= \frac{C(N,s)}{4}\|u\|^2 + I_\lambda(u) >0,
	\end{equation*}
	we have proved the claim \ref{eq4.3}.
	Next we infer that there exist $\lambda_0, t_0 > 0$ such that $J^{+}_\lambda (t\phi_{1,s}) \leq 0$ for all $t \geq t_0$ and $\lambda \in (0, \lambda_0)$.
	\begin{equation*}
		\begin{aligned}
			J^{+}_\lambda (t\phi_{1,s}) &= \frac{t^{2}}{2}\lambda_{1,s}\int\limits_{\Omega}\phi_{1,s}^{2}\,dx + \frac{\lambda}{q}t^{q}\int\limits_{\Omega}\phi_{1,s}^{q}\,dx - \frac{a}{2}t^{2}\int\limits_{\Omega}\phi_{1,s}^{2}\,dx - \frac{b}{2\cdot{2^*_{\mu ,s}}}t^{2\cdot{2^*_{\mu ,s} }}\|\phi_{1,s}\|_0^{2\cdot{2^*_{\mu ,s}}}\\
			& =\frac{t^{2}}{2}\left[\lambda_{1,s}-a \right]+ \frac{\lambda}{q}t^{q}\int\limits_{\Omega}\phi_{1,s}^{q}\,dx - \frac{b}{2\cdot{2^*_{\mu ,s}}}t^{2\cdot{2^*_{\mu ,s} }}\|\phi_{1,s}\|_0^{2\cdot{2^*_{\mu ,s}}}.
		\end{aligned}
	\end{equation*}
	As $q< 2 < 2\cdot{2^*_{\mu ,s}}$ and $a > \lambda_{1,s}$ thus we prove the claim for $\lambda>0$ small enough.
	From Lemma \ref{lemma3.1}, the functional $J^+_\lambda$ has Palais-Smale condition 
	at any level $$ c < \left( \dfrac{2^{*}_{\mu,s}-1}{2\cdot2^{*}_{\mu,s}} \right)\left[\dfrac{C(N,s)S_s^H}{2}\right]^{\frac{2^{*}_{\mu,s}}{2^{*}_{\mu,s}-1}}b^{\frac{-1}{2^{*}_{\mu,s}-1}}. $$ Set the minimax value
	$$c_\lambda := \inf\limits_{h\in \Gamma}\max\limits_{t \in [0, 1]}J^+_\lambda(h(t)),$$
	where \begin{equation*}
		\Gamma = \{h \in C\left( [0,1], X_0(\Omega)\right):\, h(0)= 0\;\text{and}\; h(1)= t_0\phi_{1,s}\}.
	\end{equation*}
Hence, \begin{equation*}\label{eq4.6}
	c_\lambda \leq \max\limits_{t \in [0, t_0]}J^+_\lambda(t\phi_{1,s})\leq \frac{\lambda t_0}{q}\int\limits_{\Omega}\phi_{1,s}^{q}\,dx.
\end{equation*}
	Choose $\lambda \in (0, \lambda_0)$, say $\lambda_1$ such that further,
	\begin{equation*}
		 c_\lambda < \left( \dfrac{2^{*}_{\mu,s}-1}{2\cdot2^{*}_{\mu,s}} \right)\left[\dfrac{C(N,s)S_s^H}{2}\right]^{\frac{2^{*}_{\mu,s}}{2^{*}_{\mu,s}-1}}b^{\frac{-1}{2^{*}_{\mu,s}-1}}.
	\end{equation*}
	Thus, for $\lambda \in (0,\lambda_1)$, the functional $J^+_\lambda$ verifies Palais-Smale condition. Moreover, summing up all the observation we see Mountain-pass Lemma holds true and it follows that $c_\lambda$ is the critical value of $J^+_\lambda$. Hence, there exists a non-trivial solution say $u_1\in X_0(\Omega)$ such that $J^+_\lambda(u_1) = c_\lambda \geq \alpha > 0$. Also it can be shown that it is a  positive solution by taking $\phi = u^-$ in \eqref{eq4}. Hence, $u_1$ is a non-trivial positive solution of problem ($P_\lambda$).\QED  
In order to obtain the negative solution we consider the following auxiliary problem,
\begin{equation}\label{eq4.8}  \begin{aligned}
		\displaystyle (-\Delta)^s u &= \lambda (u^-)^{q-1} - au^- \quad \text{in} \; \Omega,\\
		u & = 0\quad \text{in} \; \mathbb{R}^{N}\backslash\Omega.
	\end{aligned}
\end{equation}
The energy functional associated with the problem \eqref{eq4.8} is ${J^{-}_{\lambda } : X_0(\Omega) \rightarrow \mathbb{R}}$ defined as,
\begin{equation*}\label{eq4.9}
	\displaystyle J^{-}_\lambda (u)= \frac{C(N,s)}{4} \int\limits_{\mathbb{R}^{2N}} \frac{|u(x)-u(y)|^2}{|x-y|^{N+2s}}\,dxdy +\frac{\lambda}{q} \int\limits_{\Omega} (u^-(x))^{q} \,dx -\frac{a}{2}\int\limits_{\Omega} (u^-(x))^{2}\,dx\\
\end{equation*}
We observe that $J^-_\lambda$ is of $C^1$ class and its derivative is
\begin{equation}\label{eq5}
	\begin{aligned}
		\langle (J^{-})'_\lambda (u), \phi \rangle =& \frac{C(N,s)}{2} \int\limits_{\mathbb{R}^{2N}} \frac{(u(x)-u(y))(\phi(x)-\phi(y))}{|x-y|^{N+2s}}\,dxdy -\lambda \int\limits_{\Omega} (u^-(x))^{q-1}\phi(x) \,dx\\ & +a\int\limits_{\Omega} u^-(x)\phi(x)\,dx
	\end{aligned}
\end{equation}
thus the critical point of $J^-_\lambda$ are the non positive weak solution for problem ($P_\lambda$).
\begin{Theorem}\label{thm4.2}
	Let $\lambda > 0$. Then problem ($P_\lambda$) has a negative solution in $X_0(\Omega)$.
\end{Theorem}
	\proof We have to show that there exist positive constants $\alpha$ and $\rho$ such that,
	\begin{equation}\label{eq4.10}
		J^-_\lambda(u) \geq \alpha \quad \text{for all}\;\; u\in X_{0}(\Omega),\; \|u\| = \rho. 
	\end{equation}
	We construct a functional $I'_\lambda(u)$ and claim $u =0$ is strict local minimizer, \begin{equation*}
		I'_\lambda(u):= \frac{\lambda}{q} \int\limits_{\Omega} (u^-(x))^{q} \,dx -\frac{a}{2}\int\limits_{\Omega} (u^-(x))^{2}\,dx.
	\end{equation*}
	By \cite{Divya do-reg} it suffices to prove for $X_0 \cap C^0_s(\overline{\Omega})$. Now we proceed as in \cite{miyagaki} to prove the above claim \eqref{eq4.10} and also 
	proceeding in similar manner we infer that there exists $ t'_0 > 0$ such that $J^{-}_\lambda (-t\phi_{1,s}) \leq 0$ for all $t \geq t'_0$.
	Hence, we have got Mountain-pass geometry and moreover, Palais-Smale condition will also hold true by \cite{miyagaki}. 
    Thus, setting the minimax value
	$$c'_\lambda := \inf\limits_{h\in \Gamma}\max\limits_{t \in [0, 1]}J^-_\lambda(h(t)),$$
	where \begin{equation*}
		\Gamma' = \{h \in C\left( [0,1], X_0(\Omega)\right):\, h(0)= 0\;\text{and}\; h(1)= -t_0\phi_{1,s}\}.
	\end{equation*}
	We get, \begin{equation*}
		c'_\lambda \leq \max\limits_{t \in [0, t_0]}J^-_\lambda(-t\phi_{1,s})\leq \frac{\lambda t'_0}{q}\int\limits_{\Omega}\phi_{1,s}^{q}\,dx.
	\end{equation*}
	Hence, there exists a non-trivial solution say $u_2\in X_0(\Omega)$ such that $J^-_\lambda(u_2) = c'_\lambda \geq \alpha > 0$. Also it can be shown that it is a non-positive solution by taking $\phi = u^+$ in \eqref{eq5}. Hence, $u_2$ is a non-trivial non-positive solution of problem ($P_\lambda$).\QED
 \section{Linking geometry}
 In this section we will obtain the assumptions of Linking theorem for the functional $J_\lambda$, in order to get a Palais-Smale sequence at min-max critical level $c_s$. We consider the following orthogonal decomposition of the space $X_0$,
 $$X_0 = V_k\oplus W_k$$
 where \begin{equation*}
 	V_k= \spn\{\phi_{1,s}, \phi_{2,s}, \dots, \phi_{k,s}\}
 \end{equation*}
and \begin{equation*}
	W_k= \{u \in X_0(\Omega): \langle u, \phi_{j,s}\rangle =0, j=1,2,\dots,k \}.
\end{equation*}
$P_-$ is the orthogonal projection of $X_0$ onto $V_k$ and
$P_+$ is the orthogonal projection of $X_0$ onto $W_k$.
\begin{Proposition}
	\label{prop5.1}
	There exist constants $\alpha, \rho >0$ independent of $\lambda > 0$ such that $J_\lambda{(u)\geq \alpha}$ for all $u \in W_k$ with $\|u\|= \rho$.
	\end{Proposition}
	\proof By the minimizer property of the eigenvalue $\lambda_{k+1,s}$ on the infinite dimensional space $W_k$, we note that $\|u\|^{2}\geq \lambda_{k+1,s} \|u\|^2_{L^2}$, for $u \in W_k$.
	Also using Hardy-Littlewood-Sobolev inequality along with embedding theorem we get,
	\begin{equation*}
			\begin{aligned}
				J_\lambda(u)&= \frac{1}{2}\left[\dfrac{C(N, s)}{2}{\| u\|}^2-a\|u\|^2_{L^2}\right] + \dfrac{\lambda}{q}\|u\|^q_{L^q}-\dfrac{b}{2\cdot2^{*}_{\mu,s}}\|u^+\|^{2\cdot2^{*}_{\mu,s}}_{0}\\
				&\geq \frac{C(N,s)}{4}\left[ 1-\frac{a}{\lambda_{k+1,s}} \right]\|u\|^{2}-\frac{c b}{2\cdot2^*_{\mu,s}}\|u\|^{2\cdot2^*_{\mu,s}},
			\end{aligned}
	\end{equation*}
for some positive constant $c$. Thus, by taking $\|u\|=\rho$ small enough, we get our desired result.\QED
Without loss of generality, we assume $0 \in \Omega$  and fix $\delta >0$ such that $B_{4\delta} \subset \Omega$. Let $\eta \in C^{\infty}(\mathbb{R}^N)$ be a cut off function  such that
\begin{equation*}
\eta =	\begin{cases}
1 & \quad B_{\delta},\\
0 & \quad \mathbb{R}^N\backslash B_{2\delta},
\end{cases}
\end{equation*}
and for each $\epsilon > 0$, let  $u_\epsilon$ be defined as \begin{equation}\label{eq2.9}
u_\epsilon(x) = \eta(x)U_\epsilon(x)\quad for \; x\in \mathbb{R}^N.
\end{equation}
where $U_\epsilon$ is as defined in \eqref{minmimizer}. Now
we will construct a non-zero vector $e_\epsilon \in W_k$ for the Linking theorem,
\begin{equation}\label{eq5.1}
	e_\epsilon:= P_+u_\epsilon = u_\epsilon - P_-u_\epsilon \in W_k
\end{equation} which is a continuous function. For the non-triviality of $e_\epsilon$, we show that for every $K>0$ there exists $\epsilon(K)>0$ such that
$$0 \in \left\lbrace x \in \Omega: e_\epsilon(x)> K\right\rbrace$$
whenever $0< \epsilon \leq \epsilon(K)$. By \cite[p.286]{Chabrowski} we claim
\begin{equation}\label{eq5.2}
	\|P_-u_\epsilon\|_{L_\infty}\leq C\epsilon^{\frac{N-2s}{2}}
\end{equation} for all $\epsilon>0$ sufficiently small and with constant $C>0$. Using \eqref{eq2.9}, \eqref{eq5.1} and \eqref{eq5.2} we get
\begin{equation}\label{eq5.3}
	e_\epsilon(0)\geq \frac{1}{\|\tilde{u}\|_{L^{2^*_s}}}\epsilon^{-\frac{N-2s}{2}}-C\epsilon^{\frac{N-2s}{2}}\to \infty
\end{equation} as $\epsilon \to 0$. Thus, there exists $\epsilon_0>0$ such that 
$$e_\epsilon \neq 0\quad \text{for all}\;\; 0< \epsilon \leq \epsilon_0.$$
For $\epsilon \in (0, \epsilon_0]$ and $R_1,R_2 > 0$, we introduce the set
\begin{equation}\label{eq5.4}
	Q_{\epsilon, R_1, R_2}:= \left\lbrace u \in X_0(\Omega) : u= u_1+re_\epsilon, u_1 \in V_k\cap\overline{B_{R_1}}(0), 0 \leq r \leq R_2\right\rbrace. 
\end{equation} 
 $\partial Q_{\epsilon, R_1, R_2} $ denotes the relative boundary of $Q_{\epsilon, R_1, R_2}$ in the underlying finite dimensional space $V_k\oplus\mathbb{R}e_\epsilon$.

\begin{Proposition}\label{prop5.2}
There exist $R_1 > 0$ and $R_2 > 0$ sufficiently large such that
$$J_\lambda(u) \leq \frac{\lambda}{q}\int\limits_{\Omega}|u|^q\,dx,\;\; \text{for all} \; u \in \partial Q_{\epsilon, R_1, R_2} $$
for $\epsilon > 0$ sufficiently small and $\lambda > 0$.
\end{Proposition}
\proof We divide the boundary $\partial Q_{\epsilon, R_1, R_2} $ into three parts:
$\partial Q_{\epsilon, R_1, R_2} = \Gamma_1 \cup \Gamma_2 \cup \Gamma_3$, where
\begin{equation}\label{eq5.5}
	\begin{aligned}
		\Gamma_1&= V_k\cap B_{R_1},\\
		\Gamma_2&= \left\lbrace u\in X_0(\Omega): u= u_1+re_\epsilon, u_1\in V_k, \|u_1\|=R_1, 0\leq r\leq R_2\right\rbrace, \\
		\Gamma_3&= \left\lbrace u\in X_0(\Omega): u= u_1+R_2e_\epsilon, u_1\in V_k, \|u_1\|\leq R_1\right\rbrace.
	\end{aligned} 
\end{equation}
Let $u \in \Gamma_1$, from \eqref{eq5.5} we know $u\in V_k$ and also as $a > \lambda_{k,s}$ we get
\begin{equation*}
	\begin{aligned}
		J_\lambda(u) &= \frac{1}{2}\left[\dfrac{C(N, s)}{2}{\| u\|}^2-a\|u\|^2_{L^2}\right] + \dfrac{\lambda}{q}\|u\|^q_{L^q}-\dfrac{b}{2\cdot2^{*}_{\mu,s}}\|u^+\|^{2\cdot2^{*}_{\mu,s}}_{0}\\
		&\leq \frac{C(N, s)}{4}\left[ 1-\frac{a}{\lambda_{k,s}} \right]\|u\|^{2} + \frac{\lambda}{q}\int\limits_{\Omega}|u|^q\,dx\\
		&\leq\frac{\lambda}{q}\int\limits_{\Omega}|u|^q\,dx.
	\end{aligned}
\end{equation*}
Further let $u \in \Gamma_2$. By \eqref{eq5.1} and \cite[Proposition 21]{Serv} we infer that for some constant $C>0$
\begin{equation}\label{eq5.6}
	\|e_\epsilon\|^{2}\leq \|u_\epsilon\|^{2}\leq S_s^{\frac{N}{2s}}+ C\epsilon^{N-2s}.
\end{equation}
Hence, take $\delta:= \sup\limits_{0 < \epsilon \leq \epsilon_0}\|e_\epsilon\|$. In order to have $R_2 > \rho$, we must have $R_2\|e_\epsilon\|_{X_0} > \rho$ for $\rho > 0$ as defined in Proposition \ref{prop5.1}, whenever $0 <\epsilon < \epsilon_0$. We observe that $\dfrac{\rho}{\delta} < R_2$. Set $r_0= \max\{\dfrac{\rho}{\delta}, 1\}$. Thus, we get two cases,\\
\textbf{Case I}: $0\leq r \leq r_0$\\
\begin{equation*}
 \begin{aligned}
 		J_\lambda(u) 
 	& \leq \frac{1}{2}\left[ 1-\frac{a}{\lambda_{k,s}} \right]R_1^{2} + \frac{r^2}{2}\|e_\epsilon\|^{2}+ \frac{\lambda}{q}\int\limits_{\Omega}|u|^q\,dx\\ 
 	&\leq \frac{1}{2}\left[ 1-\frac{a}{\lambda_{k,s}} \right]R_1^{2} + \frac{r_0^2}{2}\delta^{2}+ \frac{\lambda}{q}\int\limits_{\Omega}|u|^q\,dx.
 \end{aligned}
\end{equation*}
Taking $R_1>0$ large enough we get,
\begin{equation*}
	J_\lambda(u)\leq \frac{\lambda}{q}\int\limits_{\Omega}|u|^q\,dx.
\end{equation*}
\textbf{Case II}: $r > r_0$\\
Without loss of generality we suppose $R_1 \geq 1$. Choose $\gamma > 0$ such that
\begin{equation}\label{eq5.7}
	\gamma < \frac{N-2s}{2(2^*_s-1)(2\cdot2^*_{\mu,s}-1)}
\end{equation}
 and since the space $V_k$ is finite dimensional, thus for some positive constant $c_0 $ and $c_1$, we denote
\begin{equation*}
	K(R_1):=\frac{1}{r_0}\sup\left\lbrace \|u_1\|_{L_\infty}: u_1\in V_k,\, \|u_1\|=R_1\right\rbrace \in[c_0R_1,\, c_1R_1].
\end{equation*}
We define
\begin{equation}\label{eq5.8}
	\Omega_\epsilon:= \left\lbrace x\in\Omega: e_\epsilon(x)>\frac{c_1}{\epsilon^{\gamma}}\right\rbrace .
\end{equation}
From \eqref{eq5.7} and \eqref{eq5.3} we infer that $0\in \Omega_\epsilon$ provided $\epsilon\in (0, \epsilon_0]$, for some $\epsilon_0>0$.
Let us suppose that $0 < \epsilon ^{\gamma} \leq \frac{1}{R_1}$. Then from \eqref{eq5.8} we have
\begin{equation*}\label{eq5.9}
	\Omega_{\epsilon, R_1}:= \{x\in\Omega: e_\epsilon(x)> K(R_1)\}\supset\Omega_\epsilon.
\end{equation*}
Observe that \begin{equation}\label{eq5.10}
	\frac{u_1(x)}{r}+e_\epsilon(x)>0\quad\text{for all}\;\;x\in \Omega_{\epsilon, R_1}.
\end{equation}
 By \eqref{eq5.6} and \eqref{eq5.10}, for some constant $C_1 >0$ we deduce that
\begin{equation*}
J_\lambda(u) \leq \frac{1}{2}\left[ 1-\frac{a}{\lambda_{k,s}} \right]R_1^{2} + \frac{r^2}{2}(S_s^{\frac{N}{2s}} + C_1)+ \frac{\lambda}{q}\int\limits_{\Omega}|u|^q\,dx - \frac{b}{2\cdot2^*_{\mu,s}}r^{2\cdot2^*_{\mu,s}}\|\frac{u_1}{r}+e_\epsilon\|^{2\cdot2^*_{\mu,s}}_{0,\Omega_{\epsilon, R_1}}.	 
\end{equation*}
Using the algebraic inequality where $a+b>0$, $b>0$, $p\geq 1$
\begin{equation*}
	(a+b)^{p}\geq |a|^{p}+ |b|^{p}-p|a|^{p-1}|b|-p|a||b|^{p-1},
\end{equation*}
on the critical choquard term we get,
\begin{equation*}
	\begin{aligned}
		&\iint\limits_{\Omega_\epsilon\times\Omega_\epsilon} \frac{(\frac{u_1}{r}+e_\epsilon)^{2^{*}_{\mu ,s}}(y)(\frac{u_1}{r}+e_\epsilon)^{2^{*}_{\mu ,s}}(x)}{|x-y|^ \mu}\,dxdy  \geq -c\iint\limits_{\Omega_\epsilon\times\Omega_\epsilon}\frac{|\frac{u_1}{r}|^{2^{*}_{\mu ,s}}(x)|\frac{u_1}{r}|^{2^{*}_{\mu ,s}-1}(y)e_\epsilon(y)}{|x-y|^ \mu}\,dxdy\\
		& -c\iint\limits_{\Omega_\epsilon\times\Omega_\epsilon}\frac{|\frac{u_1}{r}|^{2^{*}_{\mu ,s}}(x)|\frac{u_1}{r}|(y)e_\epsilon^{2^{*}_{\mu ,s}-1}(y)}{|x-y|^ \mu}\,dxdy -c\iint\limits_{\Omega_\epsilon\times\Omega_\epsilon}\frac{|e_\epsilon|^{2^{*}_{\mu ,s}}(x)|\frac{u_1}{r}|^{2^{*}_{\mu ,s}-1}(y)e_\epsilon(y)}{|x-y|^ \mu}\,dxdy\\
		&-c\iint\limits_{\Omega_\epsilon\times\Omega_\epsilon}\frac{e_\epsilon^{2^{*}_{\mu ,s}}(x)|\frac{u_1}{r}|(y)e_\epsilon^{2^{*}_{\mu ,s}-1}(y)}{|x-y|^ \mu}\,dxdy + \iint\limits_{\Omega_\epsilon\times\Omega_\epsilon} \frac{(e_\epsilon)^{2^{*}_{\mu,s}}(y)(e_\epsilon)^{2^{*}_{\mu,s}}(x)}{|x-y|^ \mu}\,dxdy.
	\end{aligned}
\end{equation*}
Using the fact that $\dfrac{u_1}{r}\in V_k$ that is $\|\dfrac{u_1}{r}\|_{L^{2^*_s}}\leq \dfrac{c}{r}\|u_1\|_{L^\infty} \leq cR_1$, applying Hardy-Littlewood-Sobolev inequality, H\"older inequality and taking the estimates $\|e_\epsilon\|_{L^{2^{*}_s-1}}$ and $\|e_\epsilon\|_{L^{2^{*}_s}}$ from \cite{ser} we get
\begin{equation*}
\begin{aligned}
	\displaystyle\iint\limits_{\Omega_\epsilon\times\Omega_\epsilon}\frac{\left|\frac{u_1}{r}\right|^{2^{*}_{\mu ,s}}(x)\left|\frac{u_1}{r}\right|^{2^{*}_{\mu ,s}-1}(y)e_\epsilon(y)}{|x-y|^ \mu}\,dxdy 
	&\leq cR_1^{2\cdot2^{*}_{\mu,s}-1}\|e_\epsilon\|_{L^{2^*_s}-1}\\
	&\leq \epsilon^{-\gamma(2\cdot2^{*}_{\mu,s}-1)}O\left( \epsilon^{\frac{N-2s}{2(2^*_s-1)}}\right).
\end{aligned}
\end{equation*}
\begin{equation*}
	\begin{aligned}
		\iint\limits_{\Omega_\epsilon\times\Omega_\epsilon}\frac{\left| \frac{u_1}{r}\right| ^{2^{*}_{\mu ,s}}(x)|\frac{u_1}{r}|(y)e_\epsilon^{2^{*}_{\mu ,s}-1}(y)}{|x-y|^ \mu}\,dxdy 
		& \leq cR_1^{2^{*}_{\mu,s}+1}{\|e_\epsilon\|^{(2^{*}_{\mu,s}-1)}_{L^{2^*_s}-1}}\\
		& \leq \epsilon^{-\gamma(2^{*}_{\mu,s}+1)}O\left( \epsilon^{\frac{(N-2s)(2^{*}_{\mu,s}-1)}{2(2^*_s-1)}}\right).
	\end{aligned}
\end{equation*}
\begin{equation*}
	\begin{aligned}
		\iint\limits_{\Omega_\epsilon\times\Omega_\epsilon}\frac{|e_\epsilon|^{2^{*}_{\mu ,s}}(x)|\frac{u_1}{r}|^{2^{*}_{\mu ,s}-1}(y)e_\epsilon(y)}{|x-y|^ \mu}\,dxdy 
		&\leq cR_1^{2^{*}_{\mu,s}-1} \|e_\epsilon\|^{2^{*}_{\mu,s}}_{L^{2^*_s}}\|e_\epsilon\|_{L^{2^*_s}-1}\\
		&\leq \left[S_s^{\frac{2N-\mu}{4s}}+O(\epsilon^{N}) \right] \epsilon^{-\gamma(2^{*}_{\mu,s}-1)}O\left( \epsilon^{\frac{N-2s}{2(2^*_s-1)}}\right).
	\end{aligned}
\end{equation*}
\begin{equation*}
	\begin{aligned}
	\iint\limits_{\Omega_\epsilon\times\Omega_\epsilon}\frac{e_\epsilon^{2^{*}_{\mu ,s}}(x)|\frac{u_1}{r}|(y)e_\epsilon^{2^{*}_{\mu ,s}-1}(y)}{|x-y|^ \mu}\,dxdy 
	& \leq cR_1\|e_\epsilon\|^{2^{*}_{\mu,s}}_{L^{2^*_s}}{\|e_\epsilon\|^{(2^{*}_{\mu,s}-1)}_{L^{2^*_s}-1}}\\
	& \leq \left[S_s^{\frac{2N-\mu}{4s}}+O(\epsilon^{N}) \right]\epsilon^{-\gamma}O\left( \epsilon^{\frac{(N-2s)(2^{*}_{\mu,s}-1)}{2(2^*_s-1)}}\right).
	\end{aligned}
\end{equation*}
By \eqref{eq5.7} we have
\begin{equation*}
 \begin{aligned}
 	\frac{(N-2s)(2^{*}_{\mu,s}-1)}{2(2^*_s-1)}-\gamma > \frac{N-2s}{2(2^*_s-1)}-\gamma(2^{*}_{\mu,s}-1) > \frac{N-2s}{2(2^*_s-1)}-\gamma(2\cdot2^{*}_{\mu,s}-1) >0,\\
 	\frac{(N-2s)(2^{*}_{\mu,s}-1)}{2(2^*_s-1)}-\gamma(2^{*}_{\mu,s}+1)  > \frac{N-2s}{2(2^*_s-1)}-\gamma(2\cdot2^{*}_{\mu,s}-1) >0.
 \end{aligned}
\end{equation*}
As $\epsilon> 0$ is sufficiently small and $2^*_{\mu,s}>2$,
 \begin{equation*}
	 \begin{aligned}
	 	\displaystyle-\frac{b}{2\cdot2^*_{\mu,s}}r^{2\cdot2^*_{\mu,s}}\lVert\frac{u_1}{r}+e_\epsilon\rVert^{2\cdot2^*_{\mu,s}}_{0,\Omega_{\epsilon}} \leq & -\frac{b}{2\cdot2^*_{\mu,s}}r^{2\cdot2^*_{\mu,s}}\left[\|e_\epsilon\|^{2\cdot2^*_{\mu,s}}_{0,\Omega_{\epsilon}} -\epsilon^{-\gamma(2\cdot2^{*}_{\mu,s}-1)}O\left( \epsilon^{\frac{N-2s}{2(2^*_s-1)}}\right)\right]\\
	 	\leq & -\frac{b}{2\cdot2^*_{\mu,s}}r^{2\cdot2^*_{\mu,s}}\left[C(N, \mu,s)(S_s)^{\frac{2N-\mu}{2s}}+ O(\epsilon^{\frac{N-2s}{2}})\right]\\  &+\frac{b}{2\cdot2^*_{\mu,s}}r^{2\cdot2^*_{\mu,s}}\left[ \epsilon^{-\gamma(2\cdot2^{*}_{\mu,s}-1)}O\left( \epsilon^{\frac{N-2s}{2(2^*_s-1)}}\right)\right].
	 \end{aligned}
\end{equation*}
Here we have used estimate of $\|u_\epsilon\|_0^{2}\,$ from \cite{Tuhina} and also the following estimate
\begin{equation*}
\left| \|u_\epsilon\|_0^{2\cdot2^{*}_{\mu,s}}-\|e_\epsilon\|_0^{2\cdot2^{*}_{\mu,s}}\right|\leq C(N,\mu,s)S_s^{\frac{2N-\mu}{2s}}+ O(\epsilon^{\frac{N-2s}{2}}),
\end{equation*}
 for some positive constant $C(N,\mu,s) >0$.
Hence, compiling all these we get
\begin{equation*}
	\begin{aligned}
J_\lambda(u) \leq& \frac{1}{2}\left[ 1-\frac{a}{\lambda_{k,s}} \right]R_1^{2} + \frac{r^2}{2}(S_s^{\frac{N}{2s}} + C_1)-\frac{b}{2\cdot2^*_{\mu,s}}r^{2\cdot2^*_{\mu,s}}\left[C(N, \mu,s)^{\frac{N}{2s}}(S_s)^{\frac{2N-\mu}{2s}} + O(\epsilon^{\frac{N-2s}{2}})\right]\\  
&+\frac{b}{2\cdot2^*_{\mu,s}}r^{2\cdot2^*_{\mu,s}}\left[ \epsilon^{-\gamma(2\cdot2^{*}_{\mu,s}-1)}O\left( \epsilon^{\frac{N-2s}{2(2^*_s-1)}}\right)\right] + \frac{\lambda}{q}\int\limits_{\Omega}|u|^q\,dx.
	\end{aligned}
\end{equation*}
Consequently, since $2\cdot2^*_{\mu, s}> 2$, there exists a positive constant $C_2 >0$ such that
\begin{equation*}
J_\lambda(u) \leq \frac{1}{2}\left[ 1-\frac{a}{\lambda_{k,s}} \right]R_1^{2} + \frac{\lambda}{q}\int\limits_{\Omega}|u|^q\,dx + C_2.
\end{equation*}
Assuming $\epsilon >0$ is arbitrarily small we have $R_1 >0$ arbitrarily large by the expression ${0< \epsilon^\gamma \leq \frac{1}{R_1}}$. Thus, we get the conclusion for $u \in \Gamma_2$.\\
Next, let $u \in \Gamma_3$, from \eqref{eq5.5} we get
\begin{equation*}
	J_\lambda(u) \leq \frac{1}{2}\left[ 1-\frac{a}{\lambda_{k,s}} \right]\|u_1\|^{2} + \frac{R_2^2}{2}\|e_\epsilon\|^{2}+ \frac{\lambda}{q}\int\limits_{\Omega}|u|^q\,dx - \frac{b}{2\cdot2^*_{\mu,s}}R_2^{2\cdot2^*_{\mu,s}}\lVert( \frac{u_1}{R_2}+e_\epsilon) ^+\rVert^{2\cdot2^*_{\mu,s}}_{0}.	 
\end{equation*}
Let us suppose that $R_1, R_2 > 0$ are arbitrarily large with $0 < \epsilon \leq \frac{R_1}{R_2} \leq 1$
\begin{equation*}
	\begin{aligned}
		J_\lambda(u) \leq & \frac{1}{2}\left[ 1-\frac{a}{\lambda_{k,s}} \right]\|u_1\|^{2} + \frac{R_2^2}{2}\left[ S_s^{\frac{N}{2s}} + C_1\left(\frac{R_1}{R_2}\right)^{N-2s}\right] -\frac{b}{2\cdot2^*_{\mu,s}}R_2^{2\cdot2^*_{\mu,s}}\lVert( \frac{u_1}{R_2}+e_\epsilon) ^+\rVert^{2\cdot2^*_{\mu,s}}_{0}\\
		& +\frac{\lambda}{q}\int\limits_{\Omega}|u|^q\,dx.
	\end{aligned}	 
\end{equation*}
As $V_k$ is finite dimensional, thus there exists constant $C_2 > 0$ such that
\begin{equation*}
	\|u_1\|_{L^\infty}\leq C_2\|u_1\|\leq C_2R_1.
\end{equation*} 
Let us define
\begin{equation*}
	\Omega_{\epsilon, R_1, R_2} := \{x \in \Omega: e_\epsilon(x) > \frac{C_2R_1}{R_2} + 1\} \supset \{x \in \Omega: e_\epsilon(x) > C_2 + 1\}:= D_\epsilon,
\end{equation*}
with the Lebesgue measure $|D_\epsilon| > 0$ provided $\epsilon > 0$ sufficiently small. On $\Omega_{\epsilon, R_1, R_2}$
\begin{equation*}
	\frac{u_1}{R_2}+ e_\epsilon > \frac{u_1}{R_2} + \frac{C_2 R_1}{R_2} + 1 \geq 1,
\end{equation*}
which implies 
\begin{equation*}
\begin{aligned}
	\frac{b}{2\cdot2^*_{\mu,s}}R_2^{2\cdot2^*_{\mu,s}}\lVert( \frac{u_1}{R_2}+e_\epsilon) ^+\rVert^{2\cdot2^*_{\mu,s}}_{0} &\geq \frac{b}{2\cdot2^*_{\mu,s}}R_2^{2\cdot2^*_{\mu,s}}\lVert( \frac{u_1}{R_2}+e_\epsilon) ^+\rVert^{2\cdot2^*_{\mu,s}}_{0, \Omega_{\epsilon, R_1, R_2}\times \Omega_{\epsilon, R_1, R_2}}\\
	&\geq\frac{b}{2\cdot2^*_{\mu,s}}R_2^{2\cdot2^*_{\mu,s}}\iint\limits_{D_\epsilon\times D_\epsilon}\frac{dxdy}{|x-y|^{\mu}}\\
	& = \frac{cb}{2\cdot2^*_{\mu,s}}R_2^{2\cdot2^*_{\mu,s}},
\end{aligned}
\end{equation*} for some constant $c > 0$. As $2\cdot2^*_{\mu, s} > 2$
\begin{equation*}
	\begin{aligned}
		J_\lambda(u) &\leq  \frac{1}{2}\left[ 1-\frac{a}{\lambda_{k,s}} \right]\|u_1\|^{2} + \frac{R_2^2}{2}\left[ S_s^{\frac{N}{2s}} + C_1\right] -\frac{cb}{2\cdot2^*_{\mu,s}}R_2^{2\cdot2^*_{\mu,s}}
		 +\frac{\lambda}{q}\int\limits_{\Omega}|u|^q\,dx\\
		 &\leq \frac{\lambda}{q}\int\limits_{\Omega}|u|^q\,dx.
	\end{aligned}	 
\end{equation*}
Last inequality holds since we have assumed $R_2 >0$ to be large enough and $R_1 >0$ such that the relation $0 < \epsilon \leq \frac{R_1}{R_2} \leq 1$ holds. \QED

\section{Third solution}

Our aim in this section is to show that the minimax value of $J_\lambda$ along $Q_{\epsilon, R_1, R_2}$ is below the admissible level given in Lemma \ref{lemma3.1}, for the compactness of  Palais-Smale sequences. Thus, by Generalized Mountain-Pass theorem we claim the third critical point of the functional, for $\epsilon >0$ and $\lambda>0$ small enough. 
\begin{Lemma}
	Let us supppose that $\epsilon> 0$ be sufficiently small and $N> 2s$, then there holds the estimate
	\begin{equation*}\label{eq5.11}
		J_\lambda(u) \leq \frac{2^*_{\mu,s}-1}{2\cdot2^*_{\mu,s}}\left(\frac{1}{b}\right)^{\frac{1}{2^*_{\mu,s}-1}} \left[\frac{\frac{C(N,s)}{2}\|e_\epsilon\|^{2}-a\|e_\epsilon\|^{2}_{L^2}}{\|u_\epsilon\|^{2}_0} \right]^{\frac{2^*_{\mu,s}}{2^*_{\mu,s}-1}}+ C\epsilon^{\frac{N-2s}{2(2^*_s-1)}}+\frac{\lambda}{q}\int\limits_{\Omega}|u|^q\,dx
	\end{equation*} for all $u \in Q_{\epsilon, R_1, R_2}$ and some positive constant $C$ independent of $\epsilon >0$.
\end{Lemma}
\proof Let $u = u_1 + re_\epsilon \in Q_{\epsilon, R_1, R_2}$, where $u_1 \in V_k\cap\overline B_{R_1}$, $0 \leq r\leq R_2$ and $e_\epsilon \in W_k$. In addition we have that $u_1$ and $e_\epsilon$ are orthogonal in $L^2$ and $X_0$ which gives,
\begin{equation}\label{5.12'}
\begin{aligned}
	J_\lambda(u) 
	 = &\frac{1}{2}\left[ \frac{C(N, s)}{2}\|u_1\|^{2}-a\|u_1\|^{2}_{L^2}\right] + \frac{r^2}{2}\left[ \frac{C(N, s)}{2}\|e_\epsilon\|^{2}-a\|e_\epsilon\|^{2}_{L^2}\right] +\frac{\lambda}{q}\int\limits_{\Omega}|u|^q\,dx\\
	& -\frac{b}{2\cdot2^*_{\mu, s}}\|u^+\|^{2\cdot2^*_{\mu, s}}_{0}\\
	\leq & \frac{r^2}{2}\left[ \frac{C(N, s)}{2}\|e_\epsilon\|^{2}-a\|e_\epsilon\|^{2}_{L^2}\right] +\frac{\lambda}{q}\int\limits_{\Omega}|u|^q\,dx - \frac{b}{2\cdot2^*_{\mu, s}}\|u^+\|^{2\cdot2^*_{\mu, s}}_{0},
\end{aligned}
\end{equation} last inequality follows as $a > \lambda_{k,s}$.
Next, we estimate the choquard term, and for that we rewrite $u$ as,
$$u= u_1 + re_\epsilon = u_1 + ru_\epsilon-rP_-u_\epsilon:= \tilde{u_1}+ ru_\epsilon.$$
Next we consider
\begin{equation*}
	\begin{aligned}
		&\int\limits_{\Omega\times\Omega}\frac{((u)^+)^{2^*_{\mu,s}}((u)^+)^{2^*_{\mu,s}}}{|x-y|^{\mu}}\,dxdy - \int\limits_{\Omega\times\Omega}\frac{(\tilde{u_1}^+)^{2^*_{\mu,s}}(\tilde{u_1}^+)^{2^*_{\mu,s}}}{|x-y|^{\mu}}\,dxdy- \int\limits_{\Omega\times\Omega}\frac{(ru_\epsilon)^{2^*_{\mu,s}}(ru_\epsilon)^{2^*_{\mu,s}}}{|x-y|^{\mu}}\,dxdy\\
		& = \int\limits_{\Omega\times\Omega}\int\limits_0^1\frac{d}{dt}\frac{((t\tilde{u_1}+ ru_\epsilon)^+)^{2^*_{\mu,s}}((t\tilde{u_1}+ ru_\epsilon)^+)^{2^*_{\mu,s}}}{|x-y|^{\mu}}-\frac{d}{dt}\frac{(t\tilde{u_1}^+)^{2^*_{\mu,s}}(t\tilde{u_1}^+)^{2^*_{\mu,s}}}{|x-y|^{\mu}}dtdxdy.
	\end{aligned}
\end{equation*}
By the symmetry of the variables and Fubini's theorem we further get
\begin{equation*}
	\begin{aligned}
		&\int\limits_{\Omega\times\Omega}\frac{((u)^+)^{2^*_{\mu,s}}((u)^+)^{2^*_{\mu,s}}}{|x-y|^{\mu}}\,dxdy - \int\limits_{\Omega\times\Omega}\frac{(\tilde{u_1}^+)^{2^*_{\mu,s}}(\tilde{u_1}^+)^{2^*_{\mu,s}}}{|x-y|^{\mu}}\,dxdy- \int\limits_{\Omega\times\Omega}\frac{(ru_\epsilon)^{2^*_{\mu,s}}(ru_\epsilon)^{2^*_{\mu,s}}}{|x-y|^{\mu}}\,dxdy\\
		&\leq 2\cdot2^*_{\mu,s}\int\limits_0^1\int\limits_{\Omega\times\Omega}\left[ \frac{((t\tilde{u_1}+ ru_\epsilon)^+(y))^{2^*_{\mu,s}}((t\tilde{u_1}+ ru_\epsilon)^+(x))^{2^*_{\mu,s}-1}}{|x-y|^{\mu}}\right]\tilde{u_1}(x)\,dxdydt\\
	    &\quad-2\cdot2^*_{\mu,s}\int\limits_0^1\int\limits_{\Omega\times\Omega}\left[\frac{(t\tilde{u_1}^+(y))^{2^*_{\mu,s}}(t\tilde{u_1}^+(x))^{2^*_{\mu,s}-1}}{|x-y|^{\mu}}\right]\tilde{u_1}(x)\,dxdydt\\
	    & = 2\cdot2^*_{\mu,s}\int\limits_0^1\int\limits_{\Omega\times\Omega}\frac{((t\tilde{u_1}+ ru_\epsilon)^+(y))^{2^*_{\mu,s}}\left[ ((t\tilde{u_1}+ ru_\epsilon)^+)^{2^*_{\mu,s}-1}-(t\tilde{u_1}^+)^{2^*_{\mu,s}-1}\right](x)\tilde{u_1}(x)}{|x-y|^{\mu}}\,dxdydt\\
	    &\quad+2\cdot2^*_{\mu,s}\int\limits_0^1\int\limits_{\Omega\times\Omega}\frac{\left[ ((t\tilde{u_1}+ ru_\epsilon)^+(y))^{2^*_{\mu,s}}-(t\tilde{u_1}^+(y))^{2^*_{\mu,s}}\right] (t\tilde{u_1}^+(x))^{2^*_{\mu,s}-1}\tilde{u_1}(x)}{|x-y|^{\mu}}\,dxdydt.
	\end{aligned}
\end{equation*}
By the mean value theorem applied to both $\psi \mapsto \psi^{2^*_{\mu,s}}$ and $\zeta \mapsto \zeta^{2^*_{\mu,s}-1}$ on $(0, \infty)$ we infer that for some $\theta_1$ and $\theta_2\in (0,1)$
\begin{equation}\label{eq5.12}
	\begin{aligned}
		&\left| \int\limits_{\Omega\times\Omega}\frac{((u)^+)^{2^*_{\mu,s}}((u)^+)^{2^*_{\mu,s}}}{|x-y|^{\mu}}\,dxdy - \int\limits_{\Omega\times\Omega}\frac{(\tilde{u_1}^+)^{2^*_{\mu,s}}(\tilde{u_1}^+)^{2^*_{\mu,s}}}{|x-y|^{\mu}}\,dxdy- \int\limits_{\Omega\times\Omega}\frac{(ru_\epsilon)^{2^*_{\mu,s}}(ru_\epsilon)^{2^*_{\mu,s}}}{|x-y|^{\mu}}\,dxdy\right| \\
		&\leq2\cdot2^*_{\mu,s}({2^*_{\mu,s}-1})\int\limits_0^1\int\limits_{\Omega\times\Omega}\frac{((t\tilde{u_1}+ ru_\epsilon)^+(y))^{2^*_{\mu,s}}\left|  t\tilde{u_1}^+(x)+\theta_1\left( (t\tilde{u_1}+ ru_\epsilon)^+-t\tilde{u_1}^+\right)(x)\right|  ^{2^*_{\mu,s}-2}}{|x-y|^{\mu}}\\
		&\quad\times\left| (t\tilde{u_1}+ ru_\epsilon)^+(x)-t\tilde{u_1}^+(x) \right||\tilde{u_1}(x)|\,dxdydt\\
		&\quad+2\cdot(2^*_{\mu,s})^{2}\int\limits_0^1\int\limits_{\Omega\times\Omega}\frac{(t\tilde{u_1}^+(x))^{2^*_{\mu,s}-1}|\tilde{u_1}(x)|\left|  t\tilde{u_1}^+(y)+\theta_2\left((t\tilde{u_1}+ ru_\epsilon)^+ -t\tilde{u_1}^+\right)(y)\right|  ^{2^*_{\mu,s}-1}}{|x-y|^{\mu}}\\
	    &\quad\times\left| (t\tilde{u_1}+ ru_\epsilon)^+(y)-t\tilde{u_1}^+(y) \right|\,dxdydt.\\
	\end{aligned}
\end{equation}
Also note that,
\begin{equation}\label{eq5.13}
	\begin{aligned}
\left|  t\tilde{u_1}^+(x)+\theta_1\left( (t\tilde{u_1}+ ru_\epsilon)^+-t\tilde{u_1}^+\right)(x)\right| &\leq \max\{(t\tilde{u_1}+ ru_\epsilon)^+(x),\, t\tilde{u_1}^+(x)\}\\
& \leq t\tilde{u_1}^+(x)+ ru_\epsilon(x)
	\end{aligned}
\end{equation}
Similarly for the other case. Next we observe that $\tau \mapsto \tau^+$ is a contraction map which gives
\begin{equation}\label{eq5.14}
\left| (t\tilde{u_1}+ ru_\epsilon)^+(x)-t\tilde{u_1}^+(x) \right|\leq ru_\epsilon(x)
\end{equation}
Substituting \eqref{eq5.13} and \eqref{eq5.14} in \eqref{eq5.12} we get
\begin{equation*}\
	\begin{aligned}
		&\left| \int\limits_{\Omega\times\Omega}\frac{((u)^+)^{2^*_{\mu,s}}((u)^+)^{2^*_{\mu,s}}}{|x-y|^{\mu}}\,dxdy - \int\limits_{\Omega\times\Omega}\frac{(\tilde{u_1}^+)^{2^*_{\mu,s}}(\tilde{u_1}^+)^{2^*_{\mu,s}}}{|x-y|^{\mu}}\,dxdy- \int\limits_{\Omega\times\Omega}\frac{(ru_\epsilon)^{2^*_{\mu,s}}(ru_\epsilon)^{2^*_{\mu,s}}}{|x-y|^{\mu}}\,dxdy\right| \\
        &\leq2\cdot2^*_{\mu,s}({2^*_{\mu,s}-1})\int\limits_0^1\int\limits_{\Omega\times\Omega}\frac{((t\tilde{u_1}+ ru_\epsilon)^+(y))^{2^*_{\mu,s}}(t\tilde{u_1}^+(x)+ ru_\epsilon(x)) ^{2^*_{\mu,s}-2}ru_\epsilon(x)|\tilde{u_1}(x)|}{|x-y|^{\mu}} \,dxdydt\\
        &\quad +2\cdot(2^*_{\mu,s})^{2}\int\limits_0^1\int\limits_{\Omega\times\Omega}\frac{(t\tilde{u_1}^+(y)+ ru_\epsilon(y))^{2^*_{\mu,s}-1}ru_\epsilon(y)(t\tilde{u_1}^+(x))^{2^*_{\mu,s}-1}|\tilde{u_1}(x)|}{|x-y|^{\mu}}\,dxdydt.\\
        \end{aligned}
        \end{equation*}
On integrating with respect to $t$ and by symmetry of the variables we have,
\begin{equation*}
\begin{aligned}
	    &\leq  c\left[ \int\limits_{\Omega\times\Omega}\frac{(\tilde{u_1}^+(y))^{2^*_{\mu,s}}(\tilde{u_1}^+(x))^{2^*_{\mu,s}-1}ru_\epsilon(x)}{|x-y|^{\mu}} +
	    \frac{(\tilde{u_1}^+(y))^{2^*_{\mu,s}}(ru_\epsilon(x))^{2^*_{\mu,s}-1}|\tilde{u_1}(x)|}{|x-y|^{\mu}}\,dxdy\right] \\
	    &\quad+ c\left[ \int\limits_{\Omega\times\Omega}\frac{(ru_\epsilon(y))^{2^*_{\mu,s}}(\tilde{u_1}^+(x))^{2^*_{\mu,s}-1}ru_\epsilon(x)}{|x-y|^{\mu}} + 
	    \frac{(ru_\epsilon(y))^{2^*_{\mu,s}}(ru_\epsilon(x))^{2^*_{\mu,s}-1}|\tilde{u_1}(x)|}{|x-y|^{\mu}}\,dxdy\right]\\
	    &\quad+ c\int\limits_{\Omega\times\Omega}\frac{(ru_\epsilon(y))^{2^*_{\mu,s}}(\tilde{u_1}^+(x))^{2^*_{\mu,s}}}{|x-y|^{\mu}}\,dxdy
       \end{aligned}
       \end{equation*}
for some positive constant c. As $u_1^+ \in V_k$, which is a finite dimensional space and using Hardy-Littlewood -Sobolev inequality and H\"older inequality we infer that
\begin{equation*}\
	\begin{aligned}
		&\left| \|u^+\|_0^{2^*_{\mu,s}} - \|\tilde{u_1}^+\|_0^{2^*_{\mu,s}} - \|ru_\epsilon\|_0^{2^*_{\mu,s}}\right|\\
		&\leq c \|\tilde{u_1}^+\|_{L^{2^*_s}}^{2^*_{\mu,s}}\|\tilde{u_1}^+\|_{L^\infty}^\frac{2^*_{\mu,s}}{2^*_s}\|\tilde{u_1}^+\|_{L^{2^*_s}-1}^{(2^*_{\mu,s}-1)-\frac{2^*_{\mu, s}}{2^*_s}}\|ru_\epsilon\|_{L^{2^*_s}-1}\\ &\quad+c\|\tilde{u_1}^+\|_{L^{2^*_s}}^{2^*_{\mu,s}} \|\tilde{u_1}^+\|_{L^\infty}^\frac{2^*_{\mu,s}}{2^*_s}\|\tilde{u_1}^+\|_{L^{2^*_s}-1}^{\frac{2^*_s-2^*_{\mu, s}}{2^*_s}}\|ru_\epsilon\|_{L^{2^*_s}-1}^{2^*_{\mu, s}-1}\\
		&\quad +c\|ru_\epsilon\|_{L^{2^*_s}}^{2^*_{\mu,s}}\|\tilde{u_1}^+\|_{L^\infty}^{\frac{2^*_{\mu,s}}{2^*_s}}\|\tilde{u_1}^+\|_{L^{2^*_s}-1}^{(2^*_{\mu,s}-1)-\frac{2^*_{\mu, s}}{2^*_s}}\|ru_\epsilon\|_{L^{2^*_s}-1}\\
		&\quad +c\|ru_\epsilon\|_{L^{2^*_s}}^{2^*_{\mu,s}}\|\tilde{u_1}^+\|_{L^\infty}^{\frac{2^*_{\mu,s}}{2^*_s}}\|\tilde{u_1}^+\|_{L^{2^*_s}-1}^{\frac{2^*_s-2^*_{\mu, s}}{2^*_s}}\|ru_\epsilon\|_{L^{2^*_s}-1}^{2^*_{\mu, s}-1} +c\|ru_\epsilon\|_{L^{2^*_s}-1}^{2^*_{\mu, s}}\|\tilde{u_1}^+\|_{L^\infty}^\frac{2\cdot2^*_{\mu,s}}{2^*_s}\|\tilde{u_1}^+\|_{L^{2^*_s}-1}^{\frac{(2^*_{s}-2)(2^*_{\mu, s})}{2^*_s}}.
\end{aligned}
\end{equation*}
Taking the estimates $\|u_\epsilon\|^{2^{*}_s-1}_{L^{2^{*}_s-1}}= O(\epsilon^{\frac{N-2s}{2}})$ and $\|u_\epsilon\|^{2^{*}_s}_{L^{2^{*}_s}} = S_s^{\frac{N}{2s}}+ O(\epsilon^N)$ from \cite{ser} we get for some positive constants $C_1, C_2, C_3, C_4, C_5$ and $C_6$,
\begin{equation*}
	\begin{aligned}
     &\left| \|u^+\|_0^{2^*_{\mu,s}} - \|\tilde{u_1}^+\|_0^{2^*_{\mu,s}} - \|ru_\epsilon\|_0^{2^*_{\mu,s}}\right|\\
     &\leq C_1\epsilon^{\frac{N-2s}{2(2^*_s-1)}} + C_2\epsilon^{\frac{(N-2s)(2^*_{\mu,s}-1)}{2(2^*_s-1)}}+ C_3[S_s^{\frac{2N-\mu}{4s}}+ O(\epsilon^N)]\epsilon^{\frac{N-2s}{2(2^*_s-1)}}\\
     &\quad +C_4[S_s^{\frac{2N-\mu}{4s}}+ O(\epsilon^N)]\epsilon^{\frac{(N-2s)(2^*_{\mu,s}-1)}{2(2^*_s-1)}}+C_5\epsilon^{\frac{(N-2s)(2^*_{\mu,s})}{2(2^*_s-1)}}\\
     & \leq C_6\epsilon^{\frac{N-2s}{2(2^*_s-1)}} .
	\end{aligned}
\end{equation*}
We have thus proved that,
\begin{equation}\label{eq5.16}
	\|u^+\|^{2\cdot2^*_{\mu,s}}_0 \geq \|\tilde{u_1}^+\|^{2\cdot2^*_{\mu,s}}_0 + \|ru_\epsilon\|^{2\cdot2^*_{\mu,s}}_0 -C_6\epsilon^{\frac{N-2s}{2(2^*_s-1)}}.
\end{equation}
From \eqref{5.12'}and \eqref{eq5.16} we conclude that
\begin{equation}\label{eq5.17}
	J_\lambda(u)\leq \frac{r^2}{2}\left[ \frac{C(N,s)}{2}\|e_\epsilon\|^{2}-a\|e_\epsilon\|^{2}_{L^2}\right] +\frac{\lambda}{q}\int\limits_{\Omega}|u|^q\,dx - \frac{b}{2\cdot2^*_{\mu, s}}r^{2\cdot2^*_{\mu, s}}\|u_\epsilon\|^{2\cdot2^*_{\mu, s}}_{0}+C_7\epsilon^{\frac{N-2s}{2(2^*_s-1)}},
\end{equation} 
for some constant $C_7>0$. Let us suppose $f: [0, \infty) \rightarrow \mathbb{R}$ is a function defined as,
\begin{equation*}
	f(t)= \frac{t^2}{2}\left[\frac{C(N,s)}{2}\|e_\epsilon\|^{2}-a\|e_\epsilon\|^{2}_{L^2} \right]- \frac{b}{2\cdot2^*_{\mu, s}}t^{2\cdot2^*_{\mu, s}}\|u_\epsilon\|^{2\cdot2^*_{\mu, s}}_{0},
\end{equation*} achieves its maximum at
$$t_0 = \left[ \frac{\frac{C(N,s)}{2}\|e_\epsilon\|^{2}-a\|e_\epsilon\|^{2}_{L^2}}{b\|u_\epsilon\|^{2\cdot2^*_{\mu,s}}_0}\right]^{\frac{1}{2\cdot2^*_{\mu, s}-2}} $$
with the maximum value
\begin{equation}\label{eq5.18}
f(t_0)= \left(\frac{2^*_{\mu, s}-1}{2\cdot2^*_{\mu, s}} \right)\left(\frac{1}{b} \right)^{\frac{1}{2^*_{\mu, s}-1}}\left[\frac{\frac{C(N,s)}{2}\|e_\epsilon\|^{2}-a\|e_\epsilon\|^{2}_{L^2}}{\|u_\epsilon\|^2_0}\right]^{\frac{2^*_{\mu, s}}{2^*_{\mu, s}-1}}.
\end{equation}
Thus, from \eqref{eq5.17} and \eqref{eq5.18} we conclude  the proof of lemma.\QED
%
%
\begin{Lemma}\label{lem6.2}
	If $N\geq 4s$, there holds the estimate 
	\begin{equation*}
		c_s:= \inf\limits_{h\in \Gamma}\max\limits_{u\in Q_{\epsilon,R_1,R_2}}J_\lambda(h(u)) <\left( \dfrac{2^{*}_{\mu,s}-1}{2\cdot2^{*}_{\mu,s}} \right)\left(\frac{1}{b}\right)^{\frac{1}{2^{*}_{\mu,s}-1}}  \left[\dfrac{C(N,s)S_s^H}{2}\right]^{\frac{2^{*}_{\mu,s}}{2^{*}_{\mu,s}-1}},
	\end{equation*}
with $\Gamma= \{h\in C(Q_{\epsilon,R_1,R_2}, X_0(\Omega)): h= \text{identity on}\;   \partial Q_{\epsilon,R_1,R_2}\}$ provided $\epsilon, \lambda >0$ are suffficiently small.
If $2s < N < 4s$, again same estimate holds, but here $k$ is sufficiently large.
\end{Lemma}
\proof Let us assume $h =$ identity on $Q_{\epsilon,R_1,R_2}$ and clearly $h \in \Gamma$. Then this implies $c_s \leq \max\limits_{u \in Q_{\epsilon,R_1,R_2}}J_\lambda(u)$. Therefore it suffices to prove that for all $u\in Q_{\epsilon,R_1,R_2}$
\begin{equation*}
	 J_\lambda(u) < \left( \dfrac{2^{*}_{\mu,s}-1}{2\cdot2^{*}_{\mu,s}} \right)\left(\frac{1}{b}\right)^{\frac{1}{2^{*}_{\mu,s}-1}}  \left[\dfrac{C(N,s)S_s^H}{2}\right]^{\frac{2^{*}_{\mu,s}}{2^{*}_{\mu,s}-1}}. 
\end{equation*}
Let us suppose $N > 4s$. By \cite{Serv}, \eqref{eq5.1} and \eqref{eq5.2} we say
\begin{equation*}
	\begin{aligned}
		\|e_\epsilon\|^{2} &\leq \|u_\epsilon\|^{2}
		\leq C(N,\mu)^{\frac{N}{2s(2^*_{\mu,s})}}(S_s^H)^{\frac{N}{2s}} +O(\epsilon^{N-2s}),\\
		\|e_\epsilon\|^{2}_{L^2} & = \|u_\epsilon\|^{2}_{L^2} - \|P_-u_\epsilon\|^{2}_{L^2}\geq C_s\epsilon^{2s}+ O(\epsilon^{N-2s})
	\end{aligned}
\end{equation*} with the constant $C_s > 0$,
\begin{equation*}
	\|u_\epsilon\|^2_0\geq \left[C(N, \mu)^{\frac{N}{2s}}(S_s^H)^{\frac{2N-\mu}{2s}}-O(\epsilon^N)\right]^{\frac{1}{2^*_{\mu,s}}}.
\end{equation*}
Next we compute for $\epsilon>0$ sufficiently small,
\begin{equation*}
	\begin{aligned}
\frac{\frac{C(N,s)}{2}\|e_\epsilon\|^{2}-a\|e_\epsilon\|^{2}_{L^2}}{\|u_\epsilon\|^2_0} &\leq \frac{\frac{C(N,s)}{2}\left[ C(N,\mu)^{\frac{N}{2s(2^*_{\mu,s})}}(S_s^H)^{\frac{N}{2s}} +O(\epsilon^{N-2s})\right]-a\left[ C_s\epsilon^{2s}+ O(\epsilon^{N-2s})\right]}{\left[C(N, \mu)^{\frac{N}{2s}}(S_s^H)^{\frac{2N-\mu}{2s}}-O(\epsilon^N)\right]^{\frac{1}{2^*_{\mu,s}}}}\\
&= \frac{C(N,s)}{2}S_s^H + O(\epsilon^{N-2s}) -aC_sO(\epsilon^{2s})\\
&< \frac{C(N,s)}{2}S_s^H.
	\end{aligned}
\end{equation*}
Next let us suppose $N= 4s$, then by \cite{Serv}
\begin{equation*}
	\|e_\epsilon\|^{2}_{L^2}  = \|u_\epsilon\|^{2}_{L^2} - \|P_-u_\epsilon\|^{2}_{L^2}\geq C_s\epsilon^{2s}|log(\epsilon)|+ O(\epsilon^{2s})
\end{equation*}
For $\epsilon, \lambda > 0$ sufficiently small, proceeding in similar manner
\begin{equation*}
	\begin{aligned}
		\frac{\frac{C(N,s)}{2}\|e_\epsilon\|^{2}-a\|e_\epsilon\|^{2}_{L^2}}{\|u_\epsilon\|^2_0} &\leq \frac{\frac{C(4s,s)}{2}\left[ C(4s,\mu)^{\frac{2}{2^*_{\mu,s}}}(S_s^H)^{2} +O(\epsilon^{2s})\right]-a\left[ C_s\epsilon^{2s}|log(\epsilon)|+ O(\epsilon^{2s})\right]}{\left[C(4s, \mu)^{2}(S_s^H)^{\frac{8s-\mu}{2s}}-O(\epsilon^{4s})\right]^{\frac{1}{2^*_{\mu,s}}}}\\
		&= \frac{C(4s,s)}{2}S_s^H + O(\epsilon^{2s}) -aC_sO(\epsilon^{2s})|log(\epsilon)|\\
		&< \frac{C(4s,s)}{2}S_s^H.
	\end{aligned}
\end{equation*}
Lastly we take the case $2s < N < 4s$, again by \cite{Serv} we have
\begin{equation*}
	\|e_\epsilon\|^{2}_{L^2}  \geq C_s\epsilon^{N-2s}+ O(\epsilon^{2s})
\end{equation*}
For $\epsilon> 0$ sufficiently small, proceeding in similar manner
\begin{equation*}
	\begin{aligned}
	\frac{\frac{C(N,s)}{2}\|e_\epsilon\|^{2}-a\|e_\epsilon\|^{2}_{L^2}}{\|u_\epsilon\|^2_0} &\leq \frac{\frac{C(N,s)}{2}\left[ C(N,\mu)^{\frac{N}{2s(2^*_{\mu,s})}}(S_s^H)^{\frac{N}{2s}} +O(\epsilon^{N-2s})\right]-a\left[ C_s\epsilon^{N-2s}+ O(\epsilon^{2s})\right]}{\left[C(N, \mu)^{\frac{N}{2s}}(S_s^H)^{\frac{2N-\mu}{2s}}-O(\epsilon^N)\right]^{\frac{1}{2^*_{\mu,s}}}}\\
	&= \frac{C(N,s)}{2}S_s^H + O(\epsilon^{N-2s}) -aC_sO(\epsilon^{N-2s})-aO(\epsilon^{2s})\\
	&< \frac{C(N,s)}{2}S_s^H.
	\end{aligned}
\end{equation*} under the assumption that $k$ is sufficiently large, where $\lambda_{k}< a < \lambda_{k+1}$. Thus the proof is complete for $\epsilon$ and $\lambda$ sufficiently small. \QED

 \textbf{Proof of Theorem \ref{thm1}}: From Theorem \ref{thm4.1} and Theorem \ref{thm4.2} we infer that we have two non-trivial solutions of opposite sign for suitable choice of $\lambda >0$, namely $u_1$ and $u_2$. Existence of third solution depends on the threshold $\overline{\lambda}$ and we apply the minimax result in Theorem \ref{thm3} to the functional $J_{\lambda} : X_0 \to \mathbb{R}$. Propostion \ref{prop5.1} verifies the condition (i) of Theorem \ref{thm3} for some $\alpha > 0$ independent of $\lambda$ and there exists some $\overline{\lambda} > 0$ such that for $\lambda \in (0, \overline{\lambda})$, Proposition \ref{prop5.2} satisfies the condition (ii) of Theorem \ref{thm3}. Since hypothesis of Theorem \ref{thm3} holds, we have a Palais-Smale sequence for the functional $J_\lambda$ at the level $c_s$ for all $\lambda \in (0, \overline{\lambda})$. Lemma \ref{lem6.2} ensures that the minimax value of $J_\lambda$ obtained along the linking
theorem is smaller than the admissible threshold for the Palais-Smale condition. Thus by Lemma \ref{lemma3.1} we conclude $J_\lambda$ satisfies the Palais-Smale condition at the level $c_s$. Hence, by Generalized Mountain-Pass theorem there exists a non-trivial critical point of the functional $J_\lambda$ say $u_3 \in X_0(\Omega)$ which is a non-trivial solution of problem $(P_\lambda) $. Lastly we show $u_3$ is distinct from $u_1$ and $u_2$.
We assume $\lambda$ is sufficiently small such that we have strict inequality in both the cases below
\begin{equation*}
	J_\lambda(u_1) = J^+_\lambda(u_1) = c_\lambda \leq \frac{\lambda t_0}{q}\int\limits_{\Omega}\phi_{1,s}^{q}\,dx < \alpha \leq c_s = J_\lambda(u_3),
\end{equation*}
	and similarly
	\begin{equation*}
		J_\lambda(u_2) = J^-_\lambda(u_2) = c'_\lambda \leq \frac{\lambda t'_0}{q}\int\limits_{\Omega}\phi_{1,s}^{q}\,dx < \alpha \leq c_s = J_\lambda(u_3),
	\end{equation*} 
then for such value of $\lambda >0$ we conclude the proof.\QED

\end{document}